\documentclass[11pt]{article}
\usepackage{latexsym}
\usepackage{epsfig}
\usepackage{amsmath}
\usepackage{amssymb}
\usepackage{amsthm}
\newtheorem{theorem}{Theorem}[section]
\newtheorem{corollary}[theorem]{Corollary}

\newtheorem{lemma}[theorem]{Lemma}

\theoremstyle{remark}
\newtheorem{remark}[theorem]{Remark}

\theoremstyle{definition}
\newtheorem{definition}[theorem]{Definition}

\DeclareMathOperator\sign{sgn}

\DeclareMathOperator\artanh{artanh}

\DeclareMathOperator\cl{cl}

\begin{document}

\title{Projectively self-concordant barriers on convex sets}

\author{Roland Hildebrand \thanks{%
Univ.\ Grenoble Alpes, CNRS, Grenoble INP, LJK, 38000 Grenoble, France
({\tt roland.hildebrand@univ-grenoble-alpes.fr}).}}

\maketitle

\begin{abstract}
Self-concordance is the most important property required for barriers in convex programming. It is intrinsically linked to the affine structure of the underlying space. Here we introduce an alternative notion of self-concordance which is linked to the projective structure. A function on a set $X \subset A^n$ in an $n$-dimensional affine space is projectively self-concordant if and only if it can be extended to an affinely self-concordant logarithmically homogeneous function on the conic extension $K \subset V^{n+1}$ of $X$ in the $(n+1)$-dimensional vector space obtained by homogenization of $A^n$. The feasible sets in conic programs, notably linear and semi-definite programs, are naturally equipped with projectively self-concordant barriers. However, the interior-point methods used to solve these programs employ only affine self-concordance. We show that estimates used in the analysis of interior-point methods are tighter for projective self-concordance, in particular inner and outer approximations of the set. This opens the way to a better tuning of parameters in interior-points algorithms to allow larger steps and hence faster convergence. Projective self-concordance is also a useful tool in the theoretical analysis of logarithmically homogeneous barriers on cones.
\end{abstract}



\section{Introduction}

A wide-spread class of algorithms used for solving convex programming problems are the interior-point methods \cite{NesNem94,Renegar01}. These methods employ barrier functions on the feasible set having a special property named \emph{self-concordance}, which was introduced by Y.E. Nesterov and A.S. Nemirovski \cite{NesNem94}.

{\definition \label{def:sc_barrier} Let $X \subset \mathbb R^n$ be a regular convex set (a closed convex set with non-empty interior and containing no lines). A \emph{self-concordant barrier} on $X$ with parameter $\nu$ is a $C^3$ function $F: X^o \to \mathbb R$ satisfying the conditions
\begin{itemize}
\item[(1)] $F''(x) \succ 0$ for all $x \in X^o$ ($X^o$ denoting the interior of $X$ and $\succ$ positive definiteness),
\item[(2)] $\lim_{x \to\partial X} F(x) = +\infty$ ($\partial X$ denoting the boundary of $X$),
\item[(3)] $|F'''(x)[h,h,h]| \leq 2(F''(x)[h,h])^{3/2}$ for all $x \in X^o$ and tangent vectors $h \in T_xX^o$,
\item[(4)] $|F'(x)[h]| \leq \sqrt{\nu}(F''(x)[h,h])^{1/2}$ for all $x \in X^o$, $h \in T_xX^o$ ($T_xX^o$ denoting the tangent space to $X^o$ at $x$).
\end{itemize}
If condition (2) is not satisfied, then we shall speak of a \emph{self-concordant function}.}

\medskip

These properties are preserved under affine transformations of $X$. Clearly the smaller the parameter $\nu$, the stronger the last condition in Definition \ref{def:sc_barrier}. Restrictions of a self-concordant barrier $F$ to intersections of $X^o$ with affine subspaces of $\mathbb R^n$ are also self-concordant barriers, with the same value of the parameter $\nu$.

In conic programming the barrier $F: K^o \to \mathbb R$ is defined on the interior of a regular convex cone $K$. Instead of the last condition in Definition \ref{def:sc_barrier} $F$ is required to satisfy a stronger condition, namely logarithmic homogeneity \cite[Prop.\ 2.3.4]{NesNem94}. A function $F$ on $K^o$ is called \emph{logarithmically homogeneous} of degree $\nu$ if
\[ F(tx) = -\nu\log t+F(x)
\]
holds for all $t > 0$ and $x \in K^o$.

The feasible set $X$ of a conic program is the intersection of the underlying regular convex cone $K$ with an affine subspace, the latter being defined by the linear equality constraints of the program. Therefore $X$ carries a self-concordant barrier with parameter $\nu$. However, this barrier has more properties than just those in Definition \ref{def:sc_barrier}. Namely, it is obtained by the restriction to $X^o$ of a logarithmically homogeneous self-concordant function with parameter $\nu$ defined on the conic extension of $X^o$. How much stronger this last condition is can be appreciated by the following result \cite{FreundJarreSchaible96}, see also \cite{CORENesterov06_30} for a slightly worse estimate with simpler proof.

\begin{theorem} \label{thm:Jarre}
Let $D \subset \mathbb R^n$ be an open convex set and $K = \{ (t,x) \in \mathbb R^{n+1} \mid t > 0,\ t^{-1}x \in D \}$ its conic extension. Let $F: D \to \mathbb R$ be a self-concordant function with parameter $\nu \geq 1$. Then the logarithmically homogeneous function
\[ F^+(t,x) = \vartheta F(t^{-1}x) - \vartheta\alpha\nu\log t
\]
with
\[ \alpha = 4+\frac{3}{\sqrt{\nu}},\ \vartheta = \frac{1}{\alpha\nu}\left( 1 + \frac{\alpha-1}{\left(\frac{\alpha}{\sqrt{\alpha}+2}-1\right)^2} \right)
\]
on $K$ is self-concordant with parameter $\nu^+ = \vartheta\alpha\nu$. The value of $\nu^+$ is the minimal possible for general $F$.
\end{theorem}

In order for condition (2) of Definition \ref{def:sc_barrier} to carry over from $F$ to $F^+$ we need also bounded-ness of the domain, because otherwise the cone $K$ has non-zero boundary points which do not correspond to boundary points of $D$.

\begin{corollary}
Let $X \subset \mathbb R^n$ be a bounded regular convex set and $K = \cl\{ (t,x) \in \mathbb R^{n+1} \mid t > 0,\ t^{-1}x \in X \}$ its conic extension. Let $F: X^o \to \mathbb R$ be a self-concordant barrier with parameter $\nu$ on $X$. Then the function
\[ F^+(t,x) = \vartheta F(t^{-1}x) - \vartheta\alpha\nu\log t
\]
with $\alpha,\vartheta,\nu^+$ as in Theorem \ref{thm:Jarre} is a self-concordant logarithmically homogeneous barrier on $K$ with parameter $\nu^+$. The value of $\nu^+$ is the minimal possible for general $F$.
\end{corollary}

The corollary states that a generic self-concordant barrier with parameter $\nu$ on $X$ can merely be extended to a logarithmically homogeneous barrier with parameter $\nu^+$, and a loss of a factor of $\upsilon(\nu) = \frac{\nu^+}{\nu}$ occurs in the parameter. The function $\upsilon(\nu)$ decreases monotonically with $\nu$ and varies in the interval $(\frac{256}{27},\frac{343+119\sqrt{7}}{27}] \approx (9.4815,24.365]$ when $\nu$ runs through $[1,+\infty)$ \cite{FreundJarreSchaible96}.

Hence the condition of extendibility to a logarithmically homogeneous self-concordant function on the conic extension of the domain of definition is different from, albeit related to self-concordance as described in Definition \ref{def:sc_barrier}. We shall call this property \emph{projective self-concordance}, to be distinguished from \emph{affine self-concordance} in Definition \ref{def:sc_barrier}.  We shall use, however, the following definition, which is independent of constructions on the conic extension.

{\definition \label{def:psc_barrier} Let $X \subset \mathbb R^n$ be a regular convex set. A \emph{projectively self-concordant barrier} on $X$ with parameter $\gamma \geq 0$ is a $C^3$ function $f: X^o \to \mathbb R$ satisfying
\begin{itemize}
\item[(1)] $f''(x) - f'(x)f'(x)^T \succ 0$ for all $x \in X^o$,
\item[(2)] $\lim_{x \to\partial X} f(x) = +\infty$,
\item[(3)] $|f'''(x)[h,h,h] - 6f''(x)[h,h]f'(x)[h] + 4(f'(x)[h])^3| \leq 2\gamma(f''(x)[h,h] - (f'(x)[h])^2)^{3/2}$ for all $x \in X^o$, $h \in T_xX^o$.
\end{itemize}
We call the trilinear symmetric form $C$ defined by the cubic polynomial in $h$ on the left-hand side of the last inequality the \emph{cubic form}, and the symmetric bilinear form $G$ defined by the quadratic polynomial in $h$ on the right-hand side the \emph{affine metric} of the barrier. }

\smallskip

Condition (3) can hence be written as $|C(x)[h,h,h]| \leq 2\gamma(G(x)[h,h])^{3/2}$. We have the following result.

{\theorem \label{thm:cone} Let $D \subset \mathbb R^n$ be an open convex set, and let $K = \{ (t,x) \in \mathbb R^{n+1} \,|\, t > 0,\ t^{-1}x \in D \}$ be its conic extension. Let $f: D \to \mathbb R$ be a $C^3$ function, let $\nu \geq 2$ be a real number, and define $F: K \to \mathbb R$ by $F(t,x) = \nu(-\log t + f(t^{-1}x))$. Set further $\gamma = \frac{\nu - 2}{\sqrt{\nu-1}}$. Then the following are equivalent:

\begin{itemize}
\item the function $f$ satisfies conditions (1) and (3) of Definition \ref{def:psc_barrier} on $D$ with parameter $\gamma$,
\item the function $F$ satisfies conditions (1) and (3) of Definition \ref{def:sc_barrier} on $K$ with parameter $\nu$.
\end{itemize} }

\smallskip

Theorem \ref{thm:cone} will be proven in Section \ref{sec:conic}. Note that condition (4) of Definition \ref{def:sc_barrier} follows from condition (1) and the logarithmic homogeneity of $F$.

Projective and affine self-concordance on a convex set are thus related as follows.

{\corollary Let $X \subset \mathbb R^n$ be a regular convex set.

If $f$ is a projectively self-concordant barrier on $X$ with parameter $\gamma$, then $F = \nu f$ is an affinely self-concordant barrier with parameter $\nu = (\gamma + \sqrt{\gamma^2+4})\sqrt{\frac{\gamma^2}{4}+1}$.

If $F$ is an affinely self-concordant barrier on $X$ with parameter $\nu$, then $f = \frac{1}{\alpha\nu}F$ is a projectively self-concordant barrier with parameter $\gamma^+ = \frac{\nu^+ - 2}{\sqrt{\nu^+-1}}$, where $\alpha,\nu^+$ depend on $\nu$ as in Theorem \ref{thm:Jarre}.}

\emph{Proof:} By Theorem \ref{thm:cone} a projectively self-concordant barrier $f$ gives rise to a logarithmically homogeneous self-concordant function on its conic extension. The affinely self-concordant barrier $F$ is then obtained by restricting this function back to $X$. The expression for $\nu$ is obtained by resolving the relation $\gamma = \frac{\nu - 2}{\sqrt{\nu-1}}$ with respect to $\nu$.

On the other hand, by Theorem \ref{thm:Jarre} an affinely self-concordant barrier $F$ gives rise to a logarithmically homogeneous self-concordant function $F^+$ on its conic extension. The projectively self-concordant barrier $f$ is then obtained by restricting the function $\frac{1}{\vartheta\alpha\nu}F^+$ to $X$.

In both cases the barrier property (2) of the generated barrier follows from the analogous property of the original barrier.
\qed

\smallskip

Let $X$ be the feasible set of a strictly feasible conic program over a cone $K$, defined by a proper affine equality constraint. If $F$ is a barrier on $K$ with parameter $\nu$, then the restriction $f$ of $\nu^{-1}F$ to $X$ is by definition projectively self-concordant with parameter $\gamma = \frac{\nu - 2}{\sqrt{\nu-1}}$. It turns out that estimates used in the analysis of interior-point methods which rely on affine self-concordance of $F$ can be improved when projective self-concordance of $f$ is taken into account. In other words, projective self-concordance incorporates in a natural way the additional information that the barrier on $X$ is obtained as a restriction of a logarithmically homogeneous barrier on the conic extension. This motivates the study of projective self-concordance on its own.


In particular, the inner approximations (Dikin ellipsoids) and outer approximations of the feasible set $X$ derived from the values of the derivatives of the barrier at the current iterate can be improved. From the derivatives of a projectively self-concordant barrier at a given point we can construct an analog of the Dikin ellipsoid, which we call \emph{Dikin set}, because albeit given by a quadratic inequality, it is not necessarily ellipsoidal. We construct also a quadratically constrained set which \emph{contains} $X$. A similar construction using affine self-concordance is known only if $X$ is bounded. Here we compare a projectively self-concordant barrier with parameter $\gamma = \frac{\nu-2}{\sqrt{\nu-1}}$ to an affinely self-concordant barrier with parameter $\nu$. These results are described in Section \ref{sec:Dikin}.

In Section \ref{sec:special} we consider bounds on the derivatives of a projectively self-concordant barrier at a point $x$ as a function of the location of $x$ with respect to the boundary of the domain of definition. We also consider a special class of barriers, namely barriers with negative curvature.

In Section \ref{sec:duality} we describe a duality theory for projectively self-concordant barriers, in which the term $\langle x,p \rangle$ in the definition of the Legendre-Fenchel dual is replaced by the term $-\log(1+\langle x,p \rangle)$. This duality theory is most complete if the domain of definition and its polar are bounded. 

In Section \ref{sec:construction} it will be shown that projectively self-concordant functions are equivariant with respect to the group of projective transformations, which is larger than the affine group corresponding to affine self-concordance. We provide means to construct projectively self-concordant barriers on more complex sets from known ones on simpler sets. 

In Section \ref{sec:examples} we consider several examples of projectively self-concordant barriers on classes of sets and individual sets.

In Section \ref{sec:Newton} we develop a short-step path-following method for projectively self-concordant barriers and compare its performance with a similar method for affinely self-concordant barriers.

Finally, in Section \ref{sec:outlook} we suggest how to apply projective self-concordance in the theoretical analysis of logarithmically homogeneous barriers on cones, and briefly consider the problem of evaluating the optimal barrier parameter of convex combinations of barriers with known parameter.



\section{Projective self-concordance and conic extension} \label{sec:conic}

Let us prove Theorem \ref{thm:cone}.

\emph{Proof} (of Theorem \ref{thm:cone}:)
We first show that the second item in the theorem implies the first one.

Let $\tilde x \in D$ and $\tilde v \in T_{\tilde x}D \setminus \{0\}$ be arbitrary, and define $x = (1,\tilde x) \in K$ and $v = (0,\tilde v) \in T_xK$. Clearly the vector $v$ is linearly independent from $x$. We have \cite[Prop. 2.3.4]{NesNem94}
\[ F'(x)[x] = -\nu, \quad F''(x)[x,x] = \nu, \quad F''[x,v] = -F'(x)[v] = -\nu f'(\tilde x)[\tilde v],
\]
and by further differentiation
\[ F'''(x)[x,x,x] = -2\nu, \quad F'''(x)[x,x,v] = -2F''(x)[x,v] = 2\nu f'(\tilde x)[\tilde v],
\]
\[ F'''(x)[x,v,v] = -2F''(x)[v,v] = -2\nu f''(\tilde x)[\tilde v,\tilde v].
\]
Consider $h = v + \alpha x \in T_xK$ for $\alpha \in \mathbb R$. We obtain
\begin{align*}
F'''(x)[h,h,h] &= \alpha^3 F'''(x)[x,x,x] + 3\alpha^2 F'''(x)[x,x,v] + 3\alpha F'''(x)[x,v,v] + F'''(x)[v,v,v] \\ &= \nu(-2\alpha^3 + 6\alpha^2f'(\tilde x)[\tilde v] - 6\alpha f''(\tilde x)[\tilde v,\tilde v] + f'''(\tilde x)[\tilde v,\tilde v,\tilde v]), \\
F''(x)[h,h] &= \alpha^2 F''(x)[x,x] + 2\alpha F''(x)[x,v] + F''(x)[v,v] = \nu(\alpha^2 - 2\alpha f'(\tilde x)[\tilde v] + f''(\tilde x)[\tilde v,\tilde v]).
\end{align*}

Now $F''(x)[h,h] > 0$ for every $\alpha$, which implies that the discriminant $(f'(\tilde x)[\tilde v])^2 - f''(\tilde x)[\tilde v,\tilde v]$ of the quadratic polynomial $\alpha^2 - 2\alpha f'(\tilde x)[\tilde v] + f''(\tilde x)[\tilde v,\tilde v]$ is negative. It follows that $f''(\tilde x) - f'(\tilde x)f'(\tilde x)^T \succ 0$. This proves the first property in Definition \ref{def:psc_barrier}.

Set $c_3 = f'''(\tilde x)[\tilde v,\tilde v,\tilde v] - 6f''(\tilde x)[\tilde v,\tilde v]f'(\tilde x)[\tilde v] + 4(f'(\tilde x)[\tilde v])^3$, $c_2 = f''(\tilde x)[\tilde v,\tilde v] - (f'(\tilde x)[\tilde v])^2 > 0$, $t = c_2^{-1/2}(\alpha - f'(\tilde x)[\tilde v])$, $\mu = c_2^{-3/2}c_3$. The self-concordance condition (3) in Definition \ref{def:sc_barrier} implies
\[ (-2\alpha^3 + 6\alpha^2f'(\tilde x)[\tilde v] - 6\alpha f''(\tilde x)[\tilde v,\tilde v] + f'''(\tilde x)[\tilde v,\tilde v,\tilde v])^2 \leq 4\nu(\alpha^2 - 2\alpha f'(\tilde x)[\tilde v] + f''(\tilde x)[\tilde v,\tilde v])^3 \qquad \forall\ \alpha \in \mathbb R
\]
or equivalently
\[ p_{\mu}(t) = 4(\nu - 1)t^6 + 12(\nu - 2)t^4 + 4\mu t^3 + 12(\nu - 3)t^2 + 12\mu t + 4\nu - \mu^2 \geq 0 \qquad \forall\ t \in \mathbb R.
\]
Let us show that this condition implies $|\mu| \leq 2\gamma$, which is the third property in Definition \ref{def:psc_barrier}. Set $\nu = \kappa^2 + 1$, $\kappa \geq 1$, then
\begin{align*}
p_{\mu}(-\kappa^{-1}) &= \left(\frac{2(\kappa^4 + 5\kappa^2 + 2)}{\kappa^3} + \mu\right)\left(\frac{2(\kappa^2-1)}{\kappa} - \mu\right), \\
p_{\mu}(\kappa^{-1}) &= \left(\frac{2(\kappa^4 + 5\kappa^2 + 2)}{\kappa^3} - \mu\right)\left(\frac{2(\kappa^2-1)}{\kappa} + \mu\right).
\end{align*}
For $\kappa \geq 1$ we have $\frac{2(\kappa^2-1)}{\kappa} \leq \frac{2(\kappa^4 + 5\kappa^2 + 2)}{\kappa^3}$, and hence both expressions $p_{\mu}(\pm\kappa^{-1})$ are simultaneously nonnegative if and only if $|\mu| \leq 2\frac{(\kappa^2-1)}{\kappa} = 2\frac{\nu-2}{\sqrt{\nu-1}} = 2\gamma$.

\medskip

Let us now show the reverse implication. By logarithmic homogeneity of $F$ we have
\[ F''(\tau x)[\tau h,\tau h] = F''(x)[h,h],\qquad F'''(\tau x)[\tau h,\tau h,\tau h] = F'''(x)[h,h,h]
\]
for all $x \in K$, $h \in T_xK$, $\tau > 0$. Therefore we need to show conditions (1) and (3) of Definition \ref{def:sc_barrier} only at $x = (1,\tilde x)$ with $\tilde x \in D$. Let $h \in T_xK \setminus \{0\}$ be arbitrary. Then there exists a unique decomposition $h = (0,\tilde v) + \alpha x$, where $\tilde v \in T_{\tilde x}D$.

Let us first consider the case $\tilde v = 0$. Then $F''(x)[h,h] = \alpha^2\nu > 0$, and $|F'''(x)[h,h,h]| = 2|\alpha|^3\nu = 2\nu^{-1/2}(F''(x)[h,h])^{3/2} < 2(F''(x)[h,h])^{3/2}$, which proves our claim.

Now consider the case $\tilde v \not= 0$. By condition (1) in Definition \ref{def:psc_barrier} the discriminant of the polynomial $\alpha^2 - 2\alpha f'(\tilde x)[\tilde v] + f''(\tilde x)[\tilde v,\tilde v]$ is negative, and hence $F''(x)[h,h] > 0$, which proves the first condition in Definition \ref{def:sc_barrier}. Assume above notations, and set $\mu_{\pm} = \pm2\gamma$. Then we get
\begin{align*}
p_{\mu_+}(t) &= 4\kappa^{-2}(\kappa t + 1)^2(\kappa(t^2 + 2)(t - 1)^2 + (\kappa - 1)((t^4 + 3t^2 + 3)\kappa + 1)), \\
p_{\mu_-}(t) &= 4\kappa^{-2}(\kappa t - 1)^2(\kappa(t^2 + 2)(t + 1)^2 + (\kappa - 1)((t^4 + 3t^2 + 3)\kappa + 1)).
\end{align*}
These polynomials are hence nonnegative by virtue of $\kappa \geq 1$. Since $p_{\mu}(t)$ is concave in $\mu$, and $\mu \in [\mu_-,\mu_+]$ by condition (3) of Definition \ref{def:psc_barrier}, the polynomial $p_{\mu}(t)$ will also be nonnegative. Reversing the chain of equivalences, we obtain the third condition in Definition \ref{def:sc_barrier}.

This completes the proof.
\qed


%
%
%

\section{Inner and outer approximations} \label{sec:Dikin}

In this section we compare the inner and outer quadratic approximations of a convex set $X$ obtained from affine and projective self-concordance of barriers $F,f$ on $X$ with parameters $\nu$ and $\gamma$, respectively, which are related by $F = \nu \cdot f$, $\gamma = \frac{\nu - 2}{\sqrt{\nu-1}}$. First we shall prove some estimates in dimension 1. The key idea to obtain these bounds is to consider the self-concordance condition as a differential inclusion giving rise to a controlled dynamical system.

{\lemma \label{lem:p_estimate} Let $I \subset \mathbb R$ be an open interval, let $p: I \to \mathbb R$ be a $C^2$ function satisfying the conditions
\begin{equation} \label{p2_equation}
p' - p^2 > 0, \qquad p'' = 6p'p - 4p^3 + 2u\gamma(p' - p^2)^{3/2},\quad u \in [-1,1]
\end{equation}
for some $\gamma \geq 0$. Here $u$ is the control, depending on the independent variable $x \in I$. Let $x_0 \in I$ be a point, and set $p_0 = p(x_0)$, $s_0 = p'(x_0)$, $g_0 = \sqrt{s_0 - p_0^2} > 0$. Let further
\[ p_{\pm}(t;p_0,s_0) = \frac{p_0 + t(g_0^2 - p_0^2 \mp \gamma g_0p_0)}{- g_0^2t^2 + (p_0t - 1)^2 \pm \gamma g_0t(p_0t - 1)},
\]
such that $p_{\pm}(x-x_0;p_0,s_0)$ are the solutions of \eqref{p2_equation} with the above initial conditions and control $u \equiv \pm1$, respectively. Let $I_{\pm} \subset \mathbb R$ be the intervals of definition of these solutions.

Then for every $x \in I \cap I_-$ we have $p_-(x-x_0;p_0,s_0) \leq p(x)$ and for every $x \in I \cap I_+$ we have $p(x) \leq p_+(x-x_0;p_0,s_0)$. }

\emph{Proof:}
It is verified by direct calculation that $p_{\pm}(0;p_0,s_0) = p_0$, $p'(0;p_0,s_0) = s_0$, and $p_{\pm}$ satisfy the above differential equation with $u \equiv \pm1$.

Let us prove the inequality $p(x) \leq p_+(x-x_0;p_0,s_0)$. At $t = 0$ the denominator of $p_+$ equals 1, but at $t_+ = \frac{1}{p_0 + \frac{\gamma}{2}g_0}$, if this value is finite, it equals $-\frac{(\gamma^2 + 4)g_0^2}{(2p_0 + \gamma g_0)^2} < 0$. Hence $x = x_0 + t_+ \not\in I_+$, and for all $x = x_0 + t \in I_+$ we have $1 - (p_0 + \frac{\gamma}{2}g_0)t > 0$. For every such $t$ we have
\[ \frac{\partial p_+}{\partial s_0} = \frac{t(1 - (p_0 + \frac{\gamma}{2}g_0)t)}{((p_0t - 1)^2 - t^2g_0^2 + \gamma t(p_0t - 1)g_0)^2},
\]
which is positive for $t > 0$ and negative for $t < 0$. The proof of the inequality is then via the Bellman principle. We have
\[ \frac{\partial p_+}{\partial p_0}s_0 + \frac{\partial p_+}{\partial s_0}(6s_0p_0 - 4p_0^3 + 2u\gamma g_0^3) - \frac{\partial p_+}{\partial t} = 2(u-1)\gamma g_0^3\frac{\partial p_+}{\partial s_0},
\]
which is non-positive for $t > 0$ and nonnegative for $t < 0$, with the extremal value zero attained only for the control $u = 1$. Therefore application of a control different from $u = 1$ only decreases the achievable maximal value of $p(x)$, and the optimal control for maximizing $p(x)$ is $u \equiv 1$.

The inequality $p_-(x-x_0;p_0,s_0) \leq p(x)$ is proven similarly by reversing the signs of $p,t,u$.
\qed

Now we easily obtain tight bounds on projectively self-concordant functions on an interval.

{\corollary \label{cor:f3_bounds} Let $I \subset \mathbb R$ be an open interval and $f: I \to \mathbb R$ a projectively self-concordant barrier on $I$ with parameter $\gamma = \frac{\nu - 2}{\sqrt{\nu - 1}}$. Let $x_0 \in I$ be a point and $f(x_0) = f_0$, $f'(x_0) = p_0$, $f''(x_0) = s_0$, $g_0 = \sqrt{s_0 - p_0^2}$. Let
\[ I_{\pm} = \left\{ x \in \mathbb R \,\left|\, (x-x_0)\left(p_0 \mp \frac{g_0}{\sqrt{\nu - 1}}\right) < 1,\ (x-x_0)\left(p_0 \pm g_0\sqrt{\nu - 1}\right) < 1 \right. \right\}
\]
be the domains of definition of the functions
\[ f_{\pm}(x) = -\frac{\nu-1}{\nu}\log\left(1-(x-x_0)\left(p_0 \mp \frac{g_0}{\sqrt{\nu - 1}}\right)\right) - \frac{1}{\nu}\log\left(1-(x-x_0)\left(p_0 \pm g_0\sqrt{\nu - 1}\right)\right) + f_0.
\]
Then for every $x \in I \cap I_-$ we have $\sigma f_-(x) \leq \sigma f(x)$ and for every $x \in I \cap I_+$ we have $\sigma f(x) \leq \sigma f_+(x)$, where $\sigma = \sign(x-x_0)$. }

\emph{Proof:}
The derivative $p = f'$ satisfies the conditions of Lemma \ref{lem:p_estimate}. It is not hard to check that
\[ f_{\pm}(x) = f_0 + \int_0^{x-x_0} p_{\pm}(t;p_0,s_0)\,dt,
\]
where $p_{\pm}$ are defined in Lemma \ref{lem:p_estimate}. Since $f(x) = f_0 + \int_{x_0}^x p(\tau)\,d\tau$, the estimates on $f$ follow from the estimates on $p$ in this lemma.
\qed

The bounds $f_{\pm}$ may escape to $+\infty$ at finite points. The domain $I$ of definition of the projectively self-concordant barrier $f$ must then contain the domain of definition of the upper bound and be contained in the domain of definition of the lower bound. This implies the following constraints on the interval $I$.

\begin{corollary} \label{cor:I_constraints}
Assume above notations.

If $p_0 + g_0\sqrt{\nu - 1} \leq 0$, or equivalently $p_0 \leq -\sqrt{\frac{\nu - 1}{\nu}s_0}$, then also $p_0 - \frac{g_0}{\sqrt{\nu - 1}} \leq 0$, and $I_+$ is unbounded to the right. It follows that $[x_0,+\infty) \subset I$. In the opposite case the right end-point of $I_+$ is given by $x_0 + \frac{1}{p_0 + g_0\sqrt{\nu - 1}}$ and hence $[x_0,x_0 + \frac{1}{p_0 + g_0\sqrt{\nu - 1}}) \subset I$.

If $p_0 + \frac{g_0}{\sqrt{\nu - 1}} \leq 0$, or equivalently $p_0 \leq -\sqrt{\frac{s_0}{\nu}}$, then also $p_0 - g_0\sqrt{\nu - 1} \leq 0$, and $I_-$ is unbounded to the right. In the opposite case the right end-point of $I_-$ is given by $x_0 + (p_0 + \frac{g_0}{\sqrt{\nu - 1}})^{-1}$, and $[x_0 + (p_0 + \frac{g_0}{\sqrt{\nu - 1}})^{-1},+\infty) \cap I = \emptyset$. \qed
\end{corollary}

\medskip

We may now proceed to the description of the inner and outer approximations of $X$.

It is well-known \cite[Theorem 2.1.1]{NesNem94} that for every interior point $x_0 \in X^o$ the \emph{Dikin ellipsoid}
\begin{equation} \label{Dikin_traditional}
E_{F,x_0} = \left\{ x \mid (x-x_0)^TF''(x_0)(x-x_0) < 1 \right\}
\end{equation}
centered on $x_0$ is an inner approximation of the set $X$. Note that it depends only on the Hessian of the barrier at the center point and hence disregards important information on the gradient of the barrier. 


We shall define the following projective analog of the Dikin ellipsoid. The \emph{Dikin set} of $f$ around the point $x_0 \in X^o$ is the set
\begin{equation} \label{Dikin_quadric}
E^p_{f,x_0} = \left\{ x \in \mathbb R^n \,\left|\, 1 - f'(x_0)[x - x_0] > \sqrt{(\nu-1)G(x_0)[x - x_0,x - x_0]}\right. \right\}.
\end{equation}
Recall that $G = f'' - f'f^{\prime T}$ is the affine metric defined in Definition \ref{def:psc_barrier}. Note that the Dikin set is not necessarily an ellipsoid, it may likewise be unbounded, with its boundary given by a paraboloid or a convex hyperboloid. The next results states that $E^p_{f,x_0}$ is an inner approximation of $X$ which at the same time contains the Dikin ellipsoid $E_{F,x_0}$.

{\lemma \label{lem:innerDikin} Let $X \subset \mathbb R^n$ be a regular convex set, $f$ a projectively self-concordant barrier on $X$ with parameter $\gamma = \frac{\nu-2}{\sqrt{\nu-1}}$, $\nu \geq 2$, and $x_0 \in X^o$ a point. Then the Dikin set $E^p_{f,x_0}$ is convex, contained in $X^o$, and contains the Dikin ellipsoid $E_{F,x_0}$, where $F = \nu \cdot f$. }

\emph{Proof:}
Let $h \not= 0$ be an arbitrary vector, and let $l$ be the line through $x_0$ parallel to $h$. Define $p_0 = \langle f'(x_0),h \rangle$, $s_0 = \langle f''(x_0)h,h \rangle$, and $g_0 = \sqrt{s_0 - p_0^2}$. Then $G(x_0)[h,h] = g_0^2$. The restriction of $f$ to the intersection $X^o \cap l$ satisfies the conditions of Corollary \ref{cor:f3_bounds}. From the consideration of the upper bound $f_+$ from this corollary we obtain the following.

If $p_0 \leq -\sqrt{\frac{\nu - 1}{\nu}s_0}$, then $X^o$ contains the whole ray $r = \{ x_0 + th \,|\, t \geq 0 \}$. If $p_0 > -\sqrt{\frac{\nu - 1}{\nu}s_0}$, then $X^o$ contains the interval $[x_0,x_0 + \frac{1}{p_0 + g_0\sqrt{\nu - 1}}h)$.


Convexity of $E^p_{f,x_0}$ follows from condition (2) in Definition \ref{def:psc_barrier}.

Let now $x = x_0 + th$, $t > 0$. The quadratic inequality defining set \eqref{Dikin_quadric} can then be rewritten as $|t^{-1}-p_0| > \sqrt{\nu-1}g_0$. Hence the intersection of the ray $r$ with the Dikin set equals the whole ray if and only if $p_0 + g_0\sqrt{\nu-1} \leq 0$. In the opposite case the ray $r$ leaves the Dikin set at the point $x_0 + (p_0 + g_0\sqrt{\nu-1})^{-1}h$. Therefore the Dikin set is contained in $X^o$.

Finally, the Dikin ellipsoid \eqref{Dikin_traditional} centered on $x_0$ is given by
\[ E_{F,x_0} = \{ x \,|\, (x-x_0)^Tf''(x_0)(x-x_0) \leq \nu^{-1} \} = \left\{ x \,\left|\, \begin{pmatrix} 1 \\ x-x_0 \end{pmatrix}^T\begin{pmatrix} 1 & 0 \\ 0 & -\nu f''(x_0) \end{pmatrix}\begin{pmatrix} 1 \\ x-x_0 \end{pmatrix}^T > 0 \right. \right\},
\]
while $E^p_{f,x_0}$ equals the connected component of $x_0$ in the set
\begin{equation} \label{Ep_intermediate}
\left\{ x \in \mathbb R^n \,\left|\, \begin{pmatrix} 1 \\ x-x_0 \end{pmatrix}^T \begin{pmatrix} 1 &\ & -f'(x_0)^T \\ -f'(x_0) &\ & \nu f'(x_0)f'(x_0)^T - (\nu-1)f''(x_0) \end{pmatrix} \begin{pmatrix} 1 \\ x-x_0 \end{pmatrix} > 0 \right. \right\}.
\end{equation}
By virtue of the matrix inequality
\[ \begin{pmatrix} 1 &\ & -f'(x_0)^T \\ -f'(x_0) &\ & \nu f'(x_0)f'(x_0)^T - (\nu-1)f''(x_0) \end{pmatrix} \succeq \frac{\nu-1}{\nu}\begin{pmatrix} 1 & 0 \\ 0 & -\nu f''(x_0) \end{pmatrix}
\]
the ${\cal S}$-lemma yields that $E_{F,x_0}$ is contained in the set \eqref{Ep_intermediate}, and therefore in $E^p_{f,x_0}$. This completes the proof.
\qed

If the set $X$ is bounded, then also an outer approximation can be constructed from the affinely self-concordant barrier $F$ on $X$ \cite[Lemma 3.2.1]{Nemirovski96Lecture}, see also \cite[Prop.~2.3.2]{NesNem94}. Namely, let $x^*$ be the minimizer of $F$ on $X$. Then the set
\[ \Gamma_{F,x^*} = \left\{ x \mid (x-x^*)^TF''(x^*)(x-x^*) \leq (\nu+2\sqrt{\nu})^2 \right\}
\]
contains $X$. The point $x^*$ is called the \emph{analytic center} of $X$. The inner and outer approximations centered at $x^*$ are hence homothetic images of each other and related by a factor of $\nu+2\sqrt{\nu}$.

{\remark Using optimal control techniques one can also construct outer approximations around arbitrary points $x_0 \in X^o$ from $F$, but these are no more defined by a quadratic inequality. }


We define the projective analog of the outer approximation of $X^o$ as
\[ \Gamma^p_{f,x_0} = \left\{ x \in \mathbb R^n \,\left|\, 1 - f'(x_0)[x-x_0] > \sqrt{\frac{1}{\nu-1}G(x_0)[x-x_0,x-x_0]} \right. \right\}.
\]
Similarly to Lemma \ref{lem:innerDikin} one proves the following result, using the lower bound $f_-$ from Corollary \ref{cor:f3_bounds}.







{\lemma \label{lem:outerDikin} Let $X \subset \mathbb R^n$ be a regular convex set, $f$ a projectively self-concordant barrier on $X$ with parameter $\gamma = \frac{\nu-2}{\sqrt{\nu-1}}$, and $x_0 \in X^o$ a point. Then the set $\Gamma^p_{f,x_0}$ defined above is convex and contains $X^o$. \qed }

\medskip

\begin{figure}
\centering
\includegraphics[width=15.03cm,height=7.01cm]{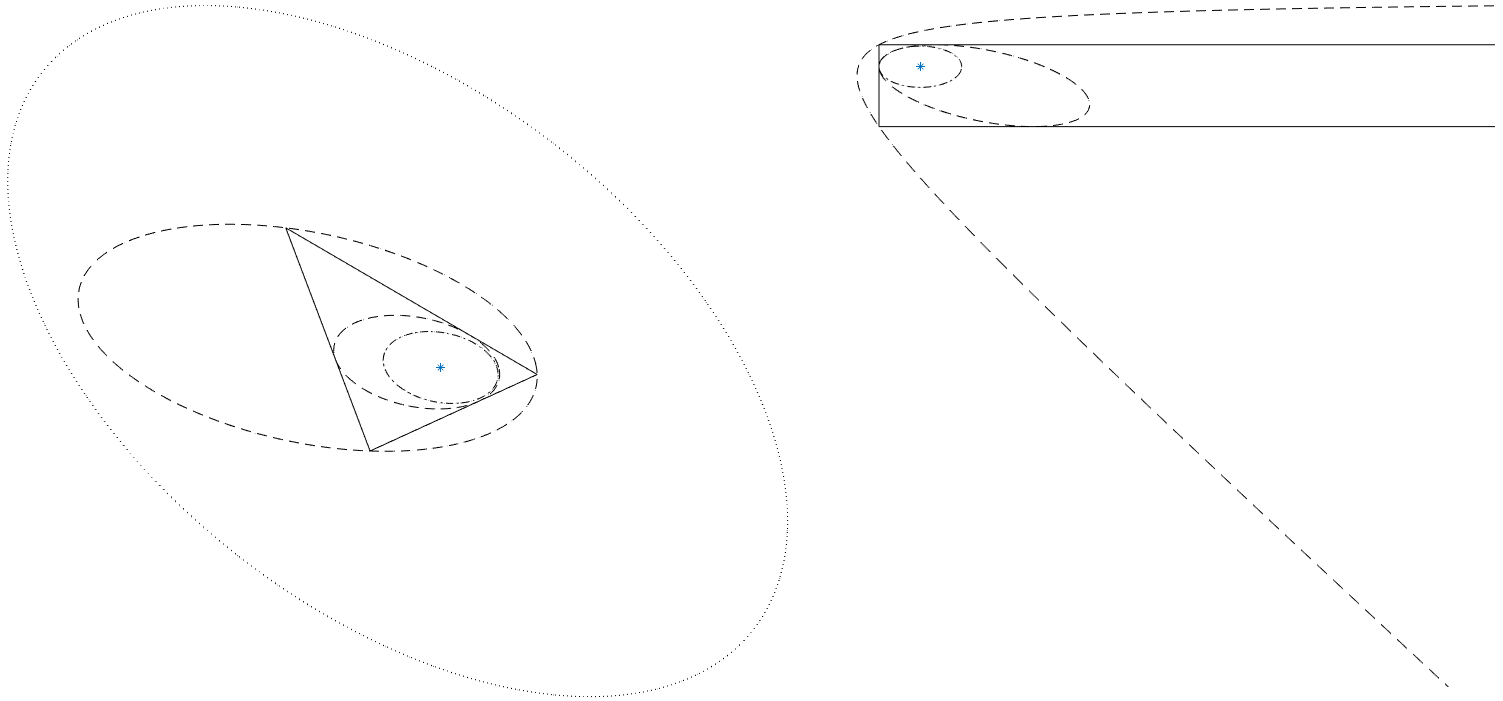}
\caption{Inner and outer approximations for a triangle (left) and a semi-infinite strip (right) around a central point $x_0$ (star): $E_{F,x_0}$ (dash-dotted), $E^p_{f,x_0}$ and $\Gamma^p_{f,x_0}$ (dashed); $\Gamma_{F,x^*}$ (dotted) around the analytic center. The barrier is the standard polyhedral logarithmic barrier, with parameter $\nu = 3$ for affine and $\gamma = \frac{\sqrt{2}}{2}$ for projective self-concordance. }
\label{fig:inner_outer_examples}
\end{figure}

In Fig.~\ref{fig:inner_outer_examples} we provide two examples how the inner and outer approximations derived from projective self-concordance of the standard logarithmic barrier for a polyhedral set relate to the inner and outer approximations computed from affine self-concordance.

Let us define the quantities
\begin{equation} \label{boundary_quantities}
\sigma_x(h) = \inf\{ t^{-1} \mid t > 0,\ x + th \in X \},\quad \pi_x(h) = \max(\sigma_x(h),\sigma_x(-h))
\end{equation}
which measure the distance to the boundary $\partial X$ of the regular convex set $X$ from a point $x \in X^o$ along a non-zero direction $h \in \mathbb R^n$. In other words, $\sigma_x,\pi_x$ are 1-homogeneous functions on $\mathbb R^n$ such that their 1-sublevel sets equal $X-x$ and $(X-x) \cap (x-X)$, accordingly. Then the inner and outer approximations determined by a projectively self-concordant function $f$ can conveniently be interpreted in terms of the quantity
\[ \omega_x(h) = \frac{\sigma_x(h) - f'(x)[h]}{\sqrt{G(x)[h,h]}}.
\]
Namely, the inclusions $E^p_{f,x_0} \subset X^o \subset \Gamma^p_{f,x_0}$ yield the following bounds on $\omega_x$.

{\corollary Let $X \subset \mathbb R^n$ be a regular convex set, and let $f$ be a projectively self-concordant barrier on $X$ with parameter $\gamma = \frac{\nu-2}{\sqrt{\nu-1}}$. Then for all $x \in X^o$ and all non-zero vectors $h \in \mathbb R^n$ we have $\frac{1}{\sqrt{\nu-1}} \leq \omega_x(h)$. If $\sigma_x(h) > 0$, then in addition $\omega_x(h) \leq \sqrt{\nu - 1}$. }

\emph{Proof:}
Let $t > 0$ be such that $t^{-1} > \sigma_x(h)$. Then $x + th \in X^o \subset \Gamma^p_{f,x_0}$, and hence $1 - f'(x)[th] > \sqrt{\frac{G(x)[th,th]}{\nu-1}}$. Dividing by $t$ and taking the limit $t^{-1} \to \sigma_x(h)$ we obtain the inequality $\omega_x(h) \geq \frac{1}{\sqrt{\nu-1}}$.

Now suppose that $\sigma_x(h) > 0$, and let $t = \frac{1}{\sigma_x(h)}$. Then $x + th \in \partial X$ and hence $x + th \not\in E^p_{f,x_0}$. It follows that $1 - f'(x)[th] \leq \sqrt{(\nu-1)G(x)[th,th]}$, which yields the desired inequality $\omega_x(h) \leq \sqrt{\nu - 1}$.
\qed

If the central point $x_0$ is the analytic center $x^*$ of $X$, then both sets $E^p_{f,x_0}$, $\Gamma^p_{f,x_0}$ are proportional to the Dikin ellipsoid $E_{F,x^*}$ and the outer approximation $\Gamma_{F,x^*}$. However, $E^p_{f,x^*}$ is $\sqrt{\frac{\nu}{\nu-1}}$ times larger and $\Gamma^p_{f,x^*}$ is $\frac{\sqrt{\nu-1}}{\sqrt{\nu}+2}$ times smaller than $E_{F,x^*}$ and $\Gamma_{F,x^*}$, respectively. 

\section{Bounds on the derivatives} \label{sec:special}

While short-step interior-point methods make steps of the order of the radius of the Dikin ellipsoid, long-step methods go a fraction of the distance to the boundary of the set. In order to justify the consistence of the latter approach, bounds on the Hessian of the barrier are needed which are valid all the way up to the boundary. Such bounds are obtained by exploiting additional structure. Here we consider such a condition involving inequalities on higher-order derivatives. First we show, however, that projective self-concordance alone also allows to derive such bounds.

In the previous section we have derived inner and outer approximations of a set carrying a projectively self-concordant barrier $f$ which are expressed in terms of the derivatives of the barrier at some point. This translates into relations between quantities \eqref{boundary_quantities} and the derivatives of $f$ at $x$.

\begin{corollary} \label{cor:bound_sigma}
Let $X \subset \mathbb R^n$ be a regular convex set, and $f$ a projectively self-concordant barrier on $X$ with parameter $\gamma = \frac{\nu-2}{\sqrt{\nu-1}}$. Then for all $x \in X^o$ and $h \in \mathbb R^n$ we have
\[ \max\left(0,f'(x)[h] + \frac{\sqrt{G(x)[h,h]}}{\sqrt{\nu-1}}\right) \leq \sigma_x(h) \leq \max\left(0,f'(x)[h] + \sqrt{\nu-1}\sqrt{G(x)[h,h]}\right),
\]
\[ |f'(x)[h]| + \frac{\sqrt{G(x)[h,h]}}{\sqrt{\nu-1}} \leq \pi_x(h) \leq |f'(x)[h]| + \sqrt{\nu-1}\sqrt{G(x)[h,h]}.
\]
\end{corollary}

\emph{Proof:}
The first chain of inequalities follows directly from Corollary \ref{cor:I_constraints}. Replacing $h$ by $-h$ we obtain
\[ \max\left(0,-f'(x)[h] + \frac{\sqrt{G(x)[h,h]}}{\sqrt{\nu-1}}\right) \leq \sigma_x(-h) \leq \max\left(0,-f'(x)[h] + \sqrt{\nu-1}\sqrt{G(x)[h,h]}\right).
\]
Combining with the first chain again, we obtain the second chain of inequalities.
\qed

In particular, we have $f'(x)[h] < \sigma_x(h)$ by the positive definiteness of the affine metric $G$. This condition can also be rewritten as
\begin{equation} \label{first_derivative_bound}
\langle f'(x),y - x \rangle < 1\quad \forall x,y \in X^o.
\end{equation}

We now consider conditions on the third and fourth derivative of a barrier $F$ involving the quantities $\sigma_x,\pi_x$.

\medskip

\emph{Barriers with negative curvature.} In \cite{NesterovTodd97} the following condition on a self-concordant barrier $F$ on a regular convex cone $K$ has been introduced:
\[ F'''(x)[h,h,u] \leq 0\quad \forall\ h \in \mathbb R^n,\ x \in K^o,\ u \in K.
\]
The condition is satisfied for the usual logarithmic barrier on spectrahedral cones and also on hyperbolicity cones \cite{Guler97}. It implies the following bound on the Hessian of the barrier. For every $x \in K^o$, $h \in \mathbb R^n$, $t \in [0,\sigma_x^{-1}(h))$ we have \cite{NesterovTodd97}
\[ \frac{1}{(1+t\sigma_x(-h))^2}F''(x) \preceq F''(x+th) \preceq \frac{1}{(1-t\sigma_x(h))^2}F''(x).
\]
Note that this condition is defined only for cones. However, this allows to define an analogous condition for projectively self-concordant barriers.

\begin{definition}
Let $X$ be a regular convex set and $f$ a projectively self-concordant barrier on $X$ with parameter $\gamma = \frac{\nu-2}{\sqrt{\nu-1}}$. We say that $f$ has \emph{negative curvature} if the $\nu$-logarithmically homogeneous extension of $\nu \cdot f$ on the conic extension $K$ of $X$ has negative curvature.
\end{definition}

Let us establish how this condition can be written in terms of the derivatives of $f$.

\begin{theorem} \label{thm:negative_projective}
A projectively self-concordant barrier $f$ on a regular convex set $X \subset \mathbb R^n$ has negative curvature if and only if for all $x \in X^o$ and $h \in \mathbb R^n$ we have
\[ \frac12(\sigma_x(h) - f'(x)[h])C(x)[\cdot,\cdot,h] \preceq (\sigma_x(h) - f'(x)[h])^2G(x) - G(x)[\cdot,h] \cdot G(x)[\cdot,h]^T.
\]
\end{theorem}

\emph{Proof:}
Let $K$ be the conic extension of $X$ and let $F$ be the logarithmically homogeneous extension of $\nu \cdot f$ to $K^o$. Here $\gamma = \frac{\nu-2}{\sqrt{\nu-1}}$ is the parameter of $f$ and $\nu$ the parameter of $F$. Let $x \in X^o$ be an arbitrary point and let $h \in T_xK^o$ be a non-zero tangent vector. Denote $\tilde x = (1,x^T)^T \in K^o$, $\tilde h = (0,h^T)^T \in T_{\tilde x}K^o$.

We have $x + th \in X$ for all $t \leq \sigma_x^{-1}(h)$. Hence $\tilde x + t\tilde h \in K$ for all $t \leq \sigma_x^{-1}(h)$ and $F'''(\tilde x)[\cdot,\cdot,\tilde x] + tF'''(\tilde x)[\cdot,\cdot,\tilde h] \preceq 0$ for all $t \leq \sigma_x^{-1}(h)$, because $F$ has negative curvature.  By logarithmic homogeneity of $F$ we have
\[ F'''(\tilde x)[\cdot,\cdot,\tilde x] = -2F''(\tilde x),\quad F''(\tilde x)[\cdot,\tilde x] = -F'(\tilde x),\quad F'(\tilde x)[\tilde x] = -\nu.
\]
It follows that $F'''(\tilde x)[\cdot,\cdot,\tilde h] \preceq 2\sigma_x(h)F''(\tilde x)$.

Let $e_1,\dots,e_n$ be a basis of the tangent space $T_xX^o$. Then $\tilde x$, $\tilde e_i = (0,e_i^T)^T$, $i = 1,\dots,n$ form a basis of the tangent space $T_{\tilde x}K^o$. In this basis we have
\[ F'''(\tilde x)[\cdot,\cdot,\tilde h] = \begin{pmatrix} 2\nu f'(x)[h] & -2\nu f''(x)[h,\cdot] \\ -2\nu f''(x)[\cdot,h] & \nu f'''(x)[\cdot,\cdot,h] \end{pmatrix},\ F''(\tilde x) = \begin{pmatrix} \nu & -\nu f'(x)^T \\ -\nu f'(x) & \nu f''(x) \end{pmatrix}.
\]
The above matrix inequality becomes
\[ \begin{pmatrix} f'(x)[h] & -f''(x)[h,\cdot] \\ -f''(x)[\cdot,h] & \frac12f'''(x)[\cdot,\cdot,h] \end{pmatrix} \preceq \sigma_x(h)\begin{pmatrix} 1 & -f'(x)^T \\ -f'(x) & f''(x) \end{pmatrix}.
\]
Since $\sigma_x(h) > f'(x)[h]$, this matrix inequality is equivalent to
\[ (\sigma_x(h)-f'(x)[h])\left(\sigma_x(h)f''(x)-\frac{f'''(x)[\cdot,\cdot,h]}{2}\right) \succeq \xi\xi^T, \quad \xi = f''(x)[\cdot,h]-\sigma_x(h)f'(x).
\]
Expressing $f''$ and $f'''$ through $G$ and $C$, respectively, we obtain the claimed relation by a somewhat lengthy calculus.
\qed

The matrix inequality in Theorem \ref{thm:negative_projective} has some interesting consequences. Applying the quadratic forms on both sides to the vector $h$ we obtain
\begin{equation} \label{omega_tangent}
C(x)[h,h,h] \leq 2\frac{\omega_x(h)^2 - 1}{\omega_x(h)}(G(x)[h,h])^{3/2}.
\end{equation}
Inequality \eqref{omega_tangent} has a number of consequences.

\begin{lemma}
Assume above notations. Then the following holds.
\begin{itemize}
\item The quantity $\omega_{x+th}(h)$ is non-decreasing as a function of $t$.
\item We have $\omega_x(h) \cdot \omega_x(-h) \geq 1$.
\item We have $\frac{\sigma_x(h) + \sigma_x(-h)}{2} \geq \sqrt{G(x)[h,h]}$.
\end{itemize}
\end{lemma}

\emph{Proof:}
Let us compute the derivative of $\omega_{x+th}(h)$ with respect to $t$. We have $\sigma_{x+th}(h) = \frac{1}{\frac{1}{\sigma_x(h)}-t}$ and hence $\frac{d\sigma_{x+th}(h)}{dt} = \sigma_{x+th}(h)^2$. Using
\begin{eqnarray*}
\frac{f'(x+th)[h]}{dt} &=& f''(x+th)[h,h] = G(x+th)[h,h] - (f'(x+th)[h])^2,\\
\frac{dG(x+th)[h,h]}{dt} &=& f'''(x+th)[h,h,h] - 2f''(x+th)[h,h]f'(x+th)[h] \\
&=& C(x+th)[h,h,h] + 4G(x+th)[h,h]f'(x+th)[h]
\end{eqnarray*}
we obtain after some calculations
\[ \frac{d\omega_{x+th}(h)}{dt} = \sqrt{G(x+th)[h,h]}\left( \omega_{x+th}(h)^2 - 1 - \frac{\omega_{x+th}(h)C(x+th)[h,h,h]}{2(G(x+th)[h,h])^{3/2}} \right) \geq 0.
\]
Here the inequality holds by \eqref{omega_tangent}. This proves the first assertion.

Replacing $h$ by $-h$ in \eqref{omega_tangent}, we obtain
\[ -C(x)[h,h,h] \leq 2\frac{\omega_x(-h)^2 - 1}{\omega_x(-h)}(G(x)[h,h])^{3/2}.
\]
Combining with \eqref{omega_tangent}, we get
\[ \frac{\omega_x(h)^2 - 1}{\omega_x(h)} + \frac{\omega_x(-h)^2 - 1}{\omega_x(-h)} \geq 0,
\]
which in view of the positivity of $\omega_x(\pm h)$ yields the second assertion.

Applying the arithmetic-geometric inequality to the second assertion we get $\frac{\omega_x(h) + \omega_x(-h)}{2} \geq 1$, which yields the third assertion.
\qed

With optimal control techniques it is possible to obtain a sharp upper bound on the Hessian $f''(x+th)$ along the ray up to the boundary. Here we provide only the bound on the value of the Hessian on the tangent direction. Denote $H_t = f''(x+th)[h,h]$, $H_0 = f''(x)[h,h]$, $\sigma = \sigma_x(h)$, $g_0 = f'(x)[h]$, then for $t \in [0,\frac{1}{\sigma})$ we have
\[ H_t \leq \frac{1}{(1 - \sigma t)^2}\left( H_0 - \frac{t(H_0 - g_0\sigma)^2[(2(1 - \sigma t) + t(\sigma - g_0))(\sigma - g_0) + t(H_0 - g_0^2)]}{(\sigma - g_0 + t(H_0 - g_0\sigma))^2} \right).
\]
It can be observed that this bound is sharper than the bound $H_t \leq \frac{H_0}{(1 - \sigma t)^2}$ obtained from the original negative curvature condition only.

%
%

\section{Duality} \label{sec:duality}

In this section we develop a duality theory for projectively self-concordant barriers. We need the following technical result.

{\lemma \label{lem:weakly_convex} Let $D \subset \mathbb R^n$ be an open convex set, and let $q: D \to \mathbb R$ be a function of class $C^2$. Suppose there exists a co-vector field $w$ on $D$ such that $q''(x) - q'(x) \cdot w(x)^T - w(x) \cdot q'(x)^T \succ 0$ for all $x \in D$. Then $q$ is quasi-convex. }

\emph{Proof:}
For the sake of contradiction, let $I \subset D$ be a closed finite interval such that $q|_I$ assumes its maximum at some point $x$ in the relative interior of $I$. Let $h \in T_xD$ be a non-zero vector in the direction of $I$. Then $q'(x)[h] = 0$, $q''(x)[h,h] \leq 0$, and hence $q''(x)[h,h] - 2q'(x)[h]\cdot w(x)[h] \leq 0$, which contradicts the assumption on $q$.

Thus for every $x,y \in D$ and every $\lambda \in (0,1)$ we have $q(\lambda x + (1-\lambda)y) \leq \max(q(x),q(y))$, i.e., $q$ is quasi-convex.
\qed

{\definition Let $X \subset \mathbb R^n$ be a regular convex set, and let $f: X^o \to \mathbb R$ be a projectively self-concordant barrier with parameter $\gamma$ on $X$. The \emph{dual} function is defined by
\[ f_*(p) = -\min_{x \in X^o}\left(f(x) + \log(1 + \langle x,p \rangle)\right).
\]
Here the domain of definition of $f_*$ consists of all points $p$ such that the function $q_p(x) = f(x) + \log(1 + \langle x,p \rangle)$ is defined and has a critical point in the interior of $X$.}

\medskip

In particular, the domain of definition of $f_*$ is contained in the interior of $-X^{\circ}$, where $X^{\circ} = \{ p \mid \langle x,p \rangle \leq 1\ \forall\ x \in X \}$ is the polar of $X$, otherwise $\log(1 + \langle x,p \rangle)$ is not defined on $X^o$.


{\lemma \label{lem:unique_min} Let $f$ be a projectively self-concordant function on a regular convex set $X$, and let $f_*$ be defined at $p$. Let $x^*$ be a critical point of the function $q_p$. Then $x^*$ is the unique critical point of $q_p$, is a global minimizer of $q_p$, and the Hessian $q_p''$ is positive definite at $x^*$. Moreover, at $x = x^*$ we have $f'(x) = -\frac{p}{1+\langle x,p \rangle}$, $p = -\frac{f'(x)}{1+\langle f'(x),x \rangle}$, $1 + \langle f'(x),x \rangle = \frac{1}{1 + \langle x,p \rangle} > 0$. }

\emph{Proof:}
The gradient of $q_p$ is given by $f'(x) + \frac{p}{1+\langle x,p \rangle}$ and the Hessian is given by
\[ q_p''(x) = f''(x) - \frac{pp^T}{(1+\langle p,x \rangle)^2} = f''(x) - f'(x)f'(x)^T + q_p'(x)w(x)^T + w(x)q_p'(x)^T
\]
with $w(x) = \frac12\left(f'(x) - \frac{p}{1+\langle x,p \rangle}\right)$. By virtue of the first property in Definition \ref{def:psc_barrier} the function $q_p$ satisfies the conditions of Lemma \ref{lem:weakly_convex} and is hence quasi-convex, and $q_p'' \succ 0$ whenever $q_p' = 0$. Every critical point of $q_p$ is hence an isolated local minimum, and there can only be one such point by quasi-convexity of $q_p$.

Now $q_p' = 0$ is equivalent to $f'(x) = -\frac{p}{1+\langle x,p \rangle}$. The other identities easily follow.
\qed

{\lemma \label{lem:duality_map} The map $x \mapsto p = -\frac{f'(x)}{1+\langle f'(x),x \rangle}$ is a bijection between the set $\{ x \in X^o \mid 1 + \langle f'(x),x \rangle > 0 \}$ and the domain of definition of $f_*$. It takes the positive definite symmetric form $f'' - f'f^{\prime T}$ to the form $f_*'' - f_*'f_*^{\prime T}$, which as a consequence is also positive definite. It also takes the symmetric 3-form $f'''_{ijk} - 2f''_{ij}f'_k - 2f''_{ik}f'_j - 2f''_{jk}f'_i + 4f'_if'_jf'_k$ to the 3-form $-(f'''_{*ijk} - 2f''_{*ij}f'_{*k} - 2f''_{*ik}f'_{*j} - 2f''_{*jk}f'_{*i} + 4f'_{*i}f'_{*j}f'_{*k})$. }

Here we denoted $f'_i = \frac{\partial f}{\partial x^i}$, $f''_{ij} = \frac{\partial^2f}{\partial x^i\partial x^j}$, $f_{*i} = \frac{\partial f_*}{\partial p_i}$ etc., $x^1,\dots,x^n$ are the coordinates of $x$, and $p_1,\dots,p_n$ the coordinates of $p$. The purpose of this notation is to be able to assume the Einstein summation convention over repeating upper and lower indices below, e.g., $\langle x,p \rangle = x^ip_i$.

\emph{Proof:}
Let $p$ be a point in the domain of definition of $f_*$. By Lemma \ref{lem:unique_min} the minimizer $x^*$ of $q_p$ is a function of $p$, namely the inverse of the map $x \mapsto p = -\frac{f'(x)}{1+\langle f'(x),x \rangle}$, and $1 + \langle f'(x^*),x^* \rangle > 0$.

On the other hand, let $x^* \in X^o$ be such that $1 + \langle f'(x^*),x^* \rangle > 0$. Set $p = -\frac{f'(x^*)}{1+\langle f'(x^*),x^* \rangle}$. Then for arbitrary $x \in X^o$ we have $1 + \langle x,p \rangle = \frac{1+\langle f'(x^*),x^* \rangle - \langle f'(x^*),x \rangle}{1+\langle f'(x^*),x^* \rangle} > 0$ by \eqref{first_derivative_bound}. Moreover, $x^*$ is a critical point of the function $q_p$ and hence $p$ is in the domain of definition of $f_*$.

This proves the first assertion of the lemma.

Differentiating the expression for $p$ we obtain
\[ \frac{\partial p}{\partial x} = \frac{(-(1+\langle f'(x),x \rangle)I + f'(x)x^T)(f''(x) - f'(x)f'(x)^T)}{(1+\langle f'(x),x \rangle)^2},
\]
and by inversion
\[ \frac{\partial x}{\partial p} = -(1+\langle f'(x),x \rangle)(f''(x) - f'(x)f'(x)^T)^{-1}(I + f'(x)x^T).
\]
The expression for $p$ also yields $f_*(p) = -f(x)+\log(1+\langle f'(x),x \rangle)$, and by differentiation
\begin{equation} \label{derivativef*}
f'_*(p) = \left(\frac{\partial x}{\partial p}\right)^T\frac{(f''(x)-f'(x)f'(x)^T)x}{1+\langle f'(x),x \rangle} = -(1+\langle f'(x),x \rangle)x.
\end{equation}

Differentiating $f_*'$ further, we get
\begin{align*}
f_*''(p) &= -\left( (1+\langle f'(x),x \rangle)I + xf'(x)^T + xx^Tf''(x) \right)\frac{\partial x}{\partial p} \\
&= -\left( (1+\langle f'(x),x \rangle)(I + xf'(x)^T) + xx^T(f''(x) - f'(x)f'(x)^T) \right)\frac{\partial x}{\partial p} \\
&= (1+\langle f'(x),x \rangle)^2\left( (I + xf'(x)^T)(f''(x) - f'(x)f'(x)^T)^{-1}(I + f'(x)x^T) + xx^T \right),
\end{align*}
and therefore
\begin{align*}
f_*''(p) - f_*'(p)f_*'(p)^T &= (1+\langle f'(x),x \rangle)^2(I + xf'(x)^T)(f''(x) - f'(x)f'(x)^T)^{-1}(I + f'(x)x^T) \\
&= \left( \frac{\partial x}{\partial p} \right)^T (f''(x) - f'(x)f'(x)^T) \frac{\partial x}{\partial p}.
\end{align*}
The second claim of the lemma follows.

Further we have
\begin{align*}
\frac{\partial x^k}{\partial p_m}\frac{\partial^2p_m}{\partial x^i\partial x^j} &= \frac{\partial x^k}{\partial p_m}\frac{\partial}{\partial x^j}\frac{(-(1+\langle f',x \rangle)\delta_m^l + f'_mx^l)(f''_{li} - f'_lf'_i)}{(1+\langle f',x \rangle)^2} \\
&= G^{km}\left( \frac{- G_{mj}G_{li}x^l - G_{mi}G_{lj}x^l}{1+\langle f',x \rangle} + f'''_{mij} - G_{mi}f'_j - G_{mj}f'_i - 2G_{ij}f'_m - 2f'_if'_jf'_m \right), \\
(f'_*)^T\frac{\partial p}{\partial x} &= \frac{x^TG}{1+\langle f'(x),x \rangle}.
\end{align*}
Here we denoted $G = f'' - f'f^{\prime T}$ and used the summation convention over repeating indices. The quantities $G_{ij}$ are the elements of the second order tensor $G$, while $G^{ij}$ are the elements of its inverse. The symbol $\delta_m^l$ is the Kronecker symbol. Differentiating the identity
\[ \left(\frac{\partial p}{\partial x}\right)^T(f_*''(p) - f_*'(p)f_*'(p)^T)\frac{\partial p}{\partial x} = f''(x) - f'(x)f'(x)^T
\]
with respect to $x$ and replacing the derivatives of $p$ by the derivatives of $f$ according to the above relations, we obtain after some calculations that
\begin{align*}
\lefteqn{\frac{\partial p_r}{\partial x^i}\frac{\partial p_s}{\partial x^j}\frac{\partial p_t}{\partial x^k}\left(f'''_{*rst} - 2f''_{*rs}f'_{*t} - 2f''_{*rt}f'_{*s} - 2f''_{*st}f'_{*r} + 4f'_{*r}f'_{*s}f'_{*t}\right) + }\\ & \qquad\qquad + f'''_{ijk} - 2f''_{ij}f'_k - 2f''_{ik}f'_j - 2f''_{jk}f'_i + 4f'_if'_jf'_k = 0.
\end{align*}
This proves the last claim.
\qed

If we denote the affine metric and the cubic form of the dual function $f_*$ by $G_*$ and $C_*$, respectively, then the lemma says that the duality map $x \mapsto p = -\frac{f'(x)}{1+\langle f'(x),x \rangle}$ takes $G$ to $G_*$ and $C$ to $-C_*$. This behaviour is similar to that of the Legendre-Fenchel duality map, which takes $F''$ to $F_*''$ and $F'''$ to $-F_*'''$ \cite[p.~45]{NesNem94}. 
Lemma \ref{lem:duality_map} easily yields the following result.

{\corollary The projective self-concordance condition on $f$ implies a similar condition on the dual function $f_*$ with the same parameter $\gamma$. \qed}

\medskip

The duality theory constructed above has the draw-back that the bijection $x \leftrightarrow p$ is only between subsets of the interiors of $X$ and $-X^{\circ}$, which may even be empty. This can be remedied if a bounded-ness assumption is introduced.

Let $X \subset \mathbb R^n$ be a regular convex set containing the origin in its interior. Setting $y = 0$ in \eqref{first_derivative_bound}, we obtain that $1 + \langle f'(x),x \rangle > 0$ for all $x \in X^o$. Bounded-ness of the polar $X^o$ hence implies that the bijection $x \leftrightarrow p$ is defined on the whole interior of $X$.

The dual assertion is also true. If $X$ is compact, then the domain of definition of the dual function $f_*$ consists of the whole interior of $-X^{\circ}$. Indeed, let $p$ be an arbitrary point from the interior of $-X^{\circ}$. Then $1 + \langle x,p \rangle$ is positive on $X$ and hence bounded below and above by positive numbers. Therefore $\log(1 + \langle x,p \rangle)$ is bounded on $X$, and the function $q_p(x)$ attains its minimum on $X^o$ by the barrier property of $f$.

We then obtain the following result.

\begin{theorem}
Let $X$ be a regular convex set such that both $X$ and its polar $X^{\circ}$ are bounded. Let $f$ be a projectively self-concordant barrier on $f$. Then the dual function $f_*(p)$ is defined on the whole interior of $-X^{\circ}$ and is actually a projectively self-concordant barrier on $-X^{\circ}$. The map $x \mapsto p = -\frac{f'(x)}{1 + \langle f'(x),x \rangle}$ is a bijection between the interiors of $X$ and $-X^{\circ}$.
\end{theorem}

\emph{Proof:}
In view of the preceding results we need only to show that $f_*$ satisfies condition (2) in Definition \ref{def:psc_barrier}. By \eqref{derivativef*} and Lemma \ref{lem:unique_min} we have $f'_*(p) = -\frac{x}{1 + \langle p,x \rangle}$, and hence the gradient $f_*'$ tends to infinity at the boundary of $-X^{\circ}$. But the function $f_*$ is affinely self-concordant, and hence must tend to infinity too.
\qed

There exists a symmetry between $f$ and $f_*$ which justifies the notion of duality. We have $f'_*(p) = -\frac{x}{1 + \langle p,x \rangle}$, or equivalently $x = -\frac{f_*'(p)}{1 + \langle f_*'(p),p \rangle}$, which is similar to the expression for $p$ as a function of $x$. Therefore the bijection between the interiors of $-X^{\circ}$ and $X$ generated by $f_*$ is the inverse of the bijection generated by $f$. It is also easily verified that the functions $F,F_*$ constructed from $f,f_*$, respectively, as in Theorem \ref{thm:cone} are the Legendre duals of each other.

\section{Construction of projectively self-concordant barriers} \label{sec:construction}

In this section we show how to construct projectively self-concordant barriers on convex sets from such barriers on simpler sets.

\medskip

\emph{Affine sections:} Let $X \subset \mathbb R^n$ be a regular convex set and $f$ a projectively self-concordant barrier on $X$ with parameter $\gamma$. Let $A \subset \mathbb R^n$ be an affine subspace intersecting the interior of $X$, and define $\tilde X = X \cap A$. From Definition \ref{def:psc_barrier} it follows in a straightforward manner that $\tilde f = f|_{\tilde X^o}$ is a projectively self-concordant barrier on $\tilde X$ with parameter $\gamma$.

\medskip

\emph{Projective images:} Let $X \subset \mathbb R^n$ be a regular convex set and $f$ a projectively self-concordant barrier on $X$ with parameter $\gamma$. Let $q$ be an affine-linear function on $\mathbb R^n$ and $A: \mathbb R^n \to \mathbb R^n$ an affine-linear isomorphism, such that there is no point at which $q$ and $A$ vanish simultaneously. Define
\[ \tilde X = \left\{ \left. \frac{A(x)}{q(x)} \,\right|\, x \in X,\ q(x) > 0 \right\}
\]
and assume this set is regular. Then the function $\tilde f$ defined by
\[ \tilde f\left(\frac{A(x)}{q(x)}\right) = f(x) + \log q(x)
\]
is a projectively self-concordant barrier on $\tilde X$ with parameter $\gamma$.

Indeed, let $y^*$ be a boundary point of $\tilde X$ and $y_k \in \tilde X^o$, $k \in \mathbb N$, a sequence of points tending to $y^*$. Let $x^* \in \partial X$ and $x_k \in X^o$ be such that $\frac{A(x^*)}{q(x^*)} = y^*$ and $\frac{A(x_k)}{q(x_k)} = y_k$. Then $q(x_k) \to q(x^*) > 0$ as $k \to +\infty$, and hence $\log q(x_k) \to \log q(x^*)$. On the other hand, $f(x_k) \to +\infty$, and hence also $\tilde f(y_k) \to +\infty$. This proves the second condition in Definition \ref{def:psc_barrier}.

The other two conditions follow from the following lemma.


{\lemma \label{lem:projective} The map $x \mapsto y = \frac{A(x)}{q(x)}$ carries the affine metric $G$ of $f$ on $X^o$ to the affine metric $\tilde G$ of $\tilde f$ on $\tilde X^o$ and the cubic form $C$ of $f$ on $X^o$ to the cubic form $\tilde C$ of $\tilde f$ on $\tilde X^o$. }

The lemma can be proven by direct calculation. Instead of reproducing this here we shall rather use a geometric result which yields an interpretation of these tensors.

\emph{Proof:}
Define the domain $D = \{ (t,x) \,|\, t > 0,\ t^{-1}x \in X^o \} \subset \mathbb R^{n+1}$ and a function $F: D \to \mathbb R$ by $F(t,x) = \log t - f(t^{-1}x)$. Let $\Gamma \subset \mathbb R^{n+1}$ be the level hypersurface $F = 0$ and let $\iota: X^o \to \Gamma$ be the bijection defined by $\iota(x) = (e^{f(x)},e^{f(x)}x)$, i.e., the map taking $x \in X^o$ to the unique point in $\Gamma$ which lies on the same ray as $(1,x)$. Then by \cite[Lemma 2.3]{Hildebrand15a} the tensors $-G$ and $-C$ on $X^o$ are taken by $\iota$ to the centro-affine metric and the centro-affine cubic form of the hypersurface $\Gamma$, respectively.

Let $\tilde D,\tilde F,\tilde\Gamma,\tilde\iota$ be similar objects defined by means of the function $\tilde f$.

Consider now the linear map $L: \mathbb R^{n+1} \to \mathbb R^{n+1}$ which takes $(1,x)$ to $(q(x),A(x))$ for every $x \in \mathbb R^n$. This map is a bijection because $q$ and $A$ do not vanish simultaneously. Moreover, the domain $\tilde D = \{ (t,y) \,|\, t > 0,\ t^{-1}y \in \tilde X^o \} \subset \mathbb R^{n+1}$ is a subset of the linear image $L[D]$. Moreover, for every $y = \frac{A(x)}{q(x)} \in \tilde X^o$ and every $\lambda > 0$ we have
\[ \tilde F(\lambda q(x),\lambda A(x)) = \log\lambda + \log q(x) - \tilde f(q(x)^{-1}A(x)) = \log\lambda + \log q(x) - f(x) - \log q(x) = F(\lambda,\lambda x).
\]
Hence $\tilde F = F \circ L^{-1}$, and the level hypersurface $\tilde\Gamma$ is a subset of $L[\Gamma]$. But linear isomorphisms leave the centro-affine metric and the centro-affine cubic form invariant by construction \cite{NomizuSasaki}. Thus $\iota^{-1} \circ L^{-1} \circ \tilde\iota$ maps the tensors $-\tilde G,-\tilde C$ to $-G,-C$, respectively. But this is exactly the map which takes $y = \frac{A(x)}{q(x)}$ to $x$, which completes the proof.
\qed

It follows that the affine metric $G$ defined by a projectively self-concordant barrier $f$ on a convex set $X$ is actually projectively invariant.

\medskip

\emph{Direct products:} Let $X_i \subset \mathbb R^{n_i}$, $i = 1,2$, be regular convex sets, and let $f_i$ be projectively self-concordant barriers on these sets with parameters $\gamma_i = \frac{\nu_i-2}{\sqrt{\nu_i-1}}$, respectively. Then $f(x,y) = \frac{\nu_1f_1(x) + \nu_2f_2(y)}{\nu_1+\nu_2}$ is a projectively self-concordant barrier on the set $X = X_1 \times X_2$ with parameter $\gamma = \frac{\nu_1+\nu_2-2}{\sqrt{\nu_1+\nu_2-1}}$.


Indeed, the second condition in Definition \ref{def:psc_barrier} follows from the corresponding conditions on the barriers $f_i$. The first and third condition follow from Theorem \ref{thm:cone} and the fact that a barrier on a direct product of cones can be constructed as a sum of barriers on the individual factor cones, its parameter being the sum of the parameters of the barriers on the factor cones.


\medskip

Note that it is possible to construct affinely self-concordant barriers on affine sections and images and direct products from an affinely self-concordant barrier on the original set. However, the equivariance with respect to projective transformations is a genuine property of projective self-concordance.

\section{Examples} \label{sec:examples}

In this section we construct projectively self-concordant barriers on different sets by virtue of the following result.

{\lemma \label{lem:cone} Let $n \geq 2$, let $K \subset \mathbb R^n$ be a regular convex cone, and $F: K^o \to \mathbb R$ a logarithmically homogeneous self-concordant barrier on $K$ with parameter $\nu$. Then $\nu \geq 2$, and $f = \nu^{-1}F|_{X^o}$ is a projectively self-concordant barrier on every proper affine section $X$ of $K$ with parameter $\gamma = \frac{\nu-2}{\sqrt{\nu-1}}$. }

\emph{Proof:}
Condition (2) in Definition \ref{def:psc_barrier} for $f$ follows from condition (2) in Definition \ref{def:sc_barrier} for $F$.

Introduce a coordinate system $(t,x_1,\dots,x_{n-1})$ in $\mathbb R^n$ such that the affine section $X$ lies in the hyperplane given by $t = 1$, and assume the notations in the first part of the proof of Theorem \ref{thm:cone}. We have $p_{\mu}(-1) + p_{\mu}(1) = 64(\nu-2) - 2\mu^2 \geq 0$, and therefore $\nu \geq 2$. By Theorem \ref{thm:cone} the first and third property in Definition \ref{def:psc_barrier} follow from the corresponding conditions on $F$ in Definition \ref{def:sc_barrier}.

This completes the proof.
\qed

{\remark If the section $X$ of $K$ is not proper, i.e., contains the origin, then $f = \nu^{-1}F|_{X^o}$ still satisfies the second and third condition in Definition \ref{def:psc_barrier}, but the matrix inequality in the first condition becomes non-strict. }

Now we are in a position to construct projectively self-concordant barriers on different sets.

\medskip

\emph{Polyhedra:} Let $P = \{ x \in \mathbb R^n \,|\, Ax \leq b \}$ be a polyhedron given by $m$ linear inequalities, with $b \not= 0$ and linearly independent columns of $A$. Then $f(x) = -\frac1m\sum_{i=1}^m \log(b-Ax)_i$ is a projectively self-concordant barrier on $P$ with parameter $\gamma = \frac{m-2}{\sqrt{m-1}}$.

Indeed, $P$ can be represented as a proper affine section of the cone $\mathbb R_+^m$, on which the standard logarithmic barrier with parameter $\nu = m$ gives rise to the above function.

\medskip

\emph{Spectrahedra:} Let $S = \{ x \,|\, {\cal A}(x) \succeq 0 \}$ be a spectrahedron given by a linear matrix inequality of size $m \times m$, with ${\cal A}$ an inhomogeneous affine map. Then $f(x) = -\frac1m\log\det{\cal A}(x)$ is a projectively self-concordant barrier on $P$ with parameter $\gamma = \frac{m-2}{\sqrt{m-1}}$.


\medskip

\emph{Epigraph of exponential function:} Consider the set $X_{\exp} = \{ (x,y) \,|\, y \geq e^x \} \subset \mathbb R^2$. On this set we have the projectively self-concordant barrier
\[ f(x,y) = -\frac13(\log(\log y - x) + \log y)
\]
with parameter $\gamma = \frac{\sqrt{2}}{2}$, which comes from the barrier on the exponential cone with parameter $\nu = 3$ defined in \cite{CharesThesis}.

\medskip

\emph{Epigraph of power functions:} For $p > 1$, consider the set $X_p = \{ (x,y) \,|\, y \geq |x|^p \}$. This set can be represented as an affine slice of the power cone
\[ K_p = \{ (x,y,z) \,|\, |x| \leq y^{1/p}z^{1/q} \},
\]
where $\frac1q = 1 - \frac1p$. The canonical barrier on this cone \cite{Hildebrand14d} leads to the projectively self-concordant barrier
\[ f(x,y) = -\frac{p+1}{3p}\log y + \frac13\phi(y^{-1/p}|x|),
\]
with the function $\phi: [0,1) \to \mathbb R$ given implicitly by the relations
\[ \log t = -\frac{1}{2p}\log\left( 1 + \frac{p+1}{\rho} \right) - \frac{1}{2q}\log\left( 1 + \frac{q+1}{\rho} \right),
\]
\[ 2\phi(t) = \left( 1 + \frac1p \right)\log(\rho+p+1) + \left( 1 + \frac1q \right)\log(\rho+q+1),
\]
with $\rho$ ranging from 0 to $+\infty$. The parameter of this barrier is given by 
\[ \gamma = \frac{\max(p,q)-2}{\sqrt{(2\max(p,q)-1)(\max(p,q)+1)}}.
\]

\medskip

We have the following general existence result.

{\corollary Let $X \subset \mathbb R^n$ be a regular convex set. Then there exists a projectively self-concordant barrier with parameter $\gamma \leq \frac{n-1}{\sqrt{n}}$ on $X$. }

\emph{Proof:}
The set $K = \cl\,\{ (t,x) \,|\, t^{-1}x \in X,\ t > 0 \} \subset \mathbb R^{n+1}$ is a regular convex cone, and $X$ can be represented as a proper affine section of $K$. But on $K$ there exist logarithmically homogeneous self-concordant barriers with parameter $\nu \leq n+1$ \cite{Hildebrand14d},\cite{BubeckEldan19}. The claim now follows from Lemma \ref{lem:cone}.
\qed

\section{A short-step path-following method} \label{sec:Newton}

In this section we analyze the performance of a Newton-like step towards a target point on the central path of an optimization problem over a set $X$ carrying a projectively self-concordant barrier. A full-blown analysis is too extensive and necessitates a separate paper, and we only sketch the algorithm and perform a quantitative analysis of the full step in one dimension. Preliminary calculation indicates, however, that the results are not too much different from those in the general case. In order to motivate the algorithm we first consider the ordinary Newton step on an affinely self-concordant function.

\medskip

Let $X \subset \mathbb R^n$ be a regular convex set and $F: X^o \to \mathbb R$ an affinely self-concordant barrier on $X$ with parameter $\nu$. We minimize the linear function $\langle c,x \rangle$ on $X$, where $c \in \mathbb R_n$ is a vector in the dual space. For $\tau \in \mathbb R_+$, let $x_*(\tau)$ be the minimizer of the composite function $F_{\tau}(x) = F(x) + \tau\langle c,x \rangle$. We suppose that this minimizer exists for large enough $\tau$. The set of minimizers is the \emph{central path} of the problem, and for $\tau \to +\infty$ it converges to an optimal solution of the original problem if the latter exists. The point $x_*(\tau)$ is hence characterized by the equation
\begin{equation} \label{affine_central_path}
F'(x_*(\tau)) = -\tau c.
\end{equation}

The simplest variant of a short-step path-following method consists in alternating a Newton step towards a target point $x_*(\tau_k)$ on the central path and an increase of the parameter $\tau_k \mapsto \tau_{k+1}$ of the target point. The ability to increase the parameter $\tau$ depends on the performance of the preceding Newton step in decreasing the \emph{Newton decrement} from the previous iterate $x_{k-1}$ to the next iterate $x_k = x_{k-1} - (F''(x_{k-1}))^{-1}F'_{\tau_k}(x_{k-1})$. Here the Newton decrement at $x$ with respect to the target parameter $\tau$ is given by
\[ \rho_{\tau}(x) = \sqrt{F'_{\tau}(x)^T(F''(x))^{-1}F'_{\tau}(x)} = ||F'(x) + \tau c||_x.
\]
The decrement is hence defined by the magnitude of the gradient mismatch between the current and the target point measured in the local Hessian metric $||\cdot||_x$ at the current point $x$.

A decrease in the decrement by $\delta_k = \rho_{\tau_k}(x_{k-1}) - \rho_{\tau_k}(x_k)$ due to the Newton step allows us to increase $\tau$ such that the decrement increases again by the same amount, i.e., $\rho_{\tau_{k+1}}(x_k) - \rho_{\tau_k}(x_k) = \delta_k$. This is equivalent to an increase of $\log\tau$ by an amount proportional to $\delta_k$ \cite[Section 3.2.4]{NesNem94}. We shall not reproduce the full analysis here, but rather provide a simplified argument. Let $v = \frac{dx_*(\tau)}{d\tau}$ be the velocity vector of the central path. Differentiating \eqref{affine_central_path} with respect to $\tau$, we obtain $F''(x_*(\tau))[v] = -c$. It follows that $v = -(F''(x_*(\tau)))^{-1}c$ and the length of the velocity vector is given by $||v||_{x_*(\tau)} = \sqrt{F''(x_*(\tau))[v,v]} = \sqrt{c^T(F''(x_*(\tau)))^{-1}c}$. On the other hand, inserting $h = (F''(x_*(\tau)))^{-1}F'(x_*(\tau))$ into condition (4) of Definition \ref{def:sc_barrier} and taking into account \eqref{affine_central_path} yields that
\begin{equation} \label{affine_length}
c^T(F''(x_*(\tau)))^{-1}c = \frac{\mu\nu}{\tau^2}, \qquad ||v||_{x_*(\tau)} = \frac{\sqrt{\mu\nu}}{\tau}
\end{equation}
for some $\mu \in (0,1]$. The length element of the central path is hence $\frac{\sqrt{\mu\nu}}{\tau}d\tau = \sqrt{\mu\nu}\,d(\log\tau)$. We may therefore increase $\log\tau$ roughly by $\frac{\delta_k}{\sqrt{\mu\nu}}$, where we neglected the dependence of $\mu$ on $\tau$ and the fact that $\delta_k$ is measured in the metric $||\cdot||_{x_k}$, while the length of the central path is measured in the metric $||\cdot||_{x_*(\tau)}$.

Let us now provide a lower bound on $\delta_k$. If $\rho_{\tau_k}(x_{k-1}) \leq \overline{\lambda}$, then the full Newton step leads to a decrement $\rho_{\tau_k}(x_k) \leq \underline{\lambda} = \left( \frac{\overline{\lambda}}{1 - \overline{\lambda}} \right)^2$ \cite[Theorem 5.2.2.1]{Nesterov18book}. One would then update the parameter $\tau$ in a way such that $\rho_{\tau_k}(x_{k-1})$ equals the maximizer $\lambda^* \approx 0.2291$ of the function $\left( \frac{\lambda}{1 - \lambda} \right)^2 - \lambda$. Then $\rho_{\tau_k}(x_k)$ is upper bounded by $\lambda_* = \left( \frac{\lambda^*}{1 - \lambda^*} \right)^2 \approx 0.0883$, and $\delta_k$ lower bounded by $\approx 0.1408$. However, the above bound on the decrement after the step is not optimal, and in order to provide a fair comparison we shall briefly sketch how to arrive at an optimal bound also for the Newton step on affinely self-concordant functions.

Let ${\cal P}$ be the set of homogeneous cubic polynomials in $\mathbb R^n$ which are bounded by 1 on the unit sphere. Each such polynomial can be seen as a symmetric third order tensor $T$. Denote by ${\cal U}$ the set $\{ U = T[\cdot,\cdot,e_1] \mid T \in {\cal P} \}$, where $e_1$ is the first basis vector in $\mathbb R^n$. Note that the elements $U$ are symmetric matrices. Then condition (3) of Definition \ref{def:sc_barrier} yields
\[ F'''(x)[\cdot,\cdot,e_1] = 2\sqrt{F''(x)[e_1,e_1]}(F''(x))^{1/2}U(F''(x))^{1/2},\qquad U \in {\cal U}.
\]
We obtain an explicit description of the evolution of the Hessian $F''$ in the direction of $e_1$.

Denote the initial point of the Newton step by $x_i$ and the final point by $x_f$. Here $x_f$ is defined as the minimizer of the quadratic approximation of the function $F_{\tau}$ at $x_i$, and the decrement $\overline{\lambda}$ at $x_i$ is the distance between $x_i$ and $x_f$ in the local metric $||\cdot||_{x_i}$. Instead of jumping from $x_i$ to $x_f$, let us travel with local unit velocity along the line segment $\sigma$ joining these points. As we move along the segment, the quadratic approximation changes, and so does its minimizer $x_m$. Denote the distance from an intermediate point $x(t) \in \sigma$ to the final point $x_f$ and the minimizer $x_m(t)$ by $\Delta_f(t),\Delta_m(t)$, respectively, and the angle between the directions toward these points by $\varphi(t)$, all measured in the local norm $||\cdot||_{x(t)}$. Affine invariance of the Newton scheme allows to obtain a controlled dynamical system governing the evolution of these quantities. Let $u_{\varphi} = (\cos\varphi,\sin\varphi,0,\dots,0)^T$. A straightforward calculus yields
\[ \frac{d\Delta_f}{dt} = -1+\Delta_f \cdot e_1^TUe_1,\quad \frac{d\Delta_m}{dt} = -\cos\varphi-\Delta_m \cdot u_{\varphi}^TUu_{\varphi}, \quad \frac{d\varphi}{dt} = \frac{\sin\varphi}{\Delta_m}+\frac{\cos\varphi}{\sin\varphi}(e_1^TUe_1 - u_{\varphi}^TUu_{\varphi}),
\]
where $U \in {\cal U}$ is the control. At the initial point, i.e., for $t = 0$, we have $\Delta_f(0) = \Delta_m(0) = \overline{\lambda}$, $\varphi(0) = 0$. When we arrive at $x_f$ at some final moment $t = T$, the distance $\Delta_f(T)$ shrinks to zero, while the distance $\Delta_m(T)$ equals the value of the decrement at the final point of the step. An upper bound $\underline{\lambda}$ on this decrement is then given by the maximal value of $\Delta_m(T)$ which can be achieved along the trajectories of this controlled system with the boundary conditions listed above.

The problem can be solved by optimal control techniques, but its solution cannot be expressed in closed form. Therefore we shall consider the simplified case of a function $F$ defined on a 1-dimensional domain. Then the set of controls ${\cal U}$ equals the interval $[-1,1]$, and the angle $\varphi$ is piece-wise constant and changes its value by $\pi$ if the minimizer $x_m(t)$ meets the intermediate point $x(t)$. Allowing for negative values of $\Delta_m$, we may dispose of the angle and arrive at the simplified system
\[ \frac{d\Delta_f}{dt} = -1+\Delta_f u,\quad \frac{d\Delta_m}{dt} = -1-\Delta_m u, \quad u \in [-1,1].
\]
This problem has the following solution. The optimal strategy maximizing $|\Delta_m(T)|$ is to apply at first the control $u \equiv -1$. This leads to the minimizer $x_m(t)$ moving closer, and the current point $x(t)$ crosses the minimum of $F_{\tau}$ and overshoots. At the time instant when $x(t)$ is in the middle between $x_m(t)$ and $x_f$ the optimal control switches to $u \equiv 1$, which has the effect of increasing the local metric and hence putting an even larger distance between $x(t)$ and $x_m(t)$. The maximal value of $|\Delta_m(T)|$ achieved by this strategy is given by
\[ \underline{\lambda} = 4 - \overline{\lambda}^2 - 4\sqrt{1 - \overline{\lambda}^2}.
\]

Finally, maximizing $\overline{\lambda} - \underline{\lambda}$ with respect to $\overline{\lambda}$ yields the optimal value $\lambda^* \approx 0.4166$ of the decrement, leading to an upper bound $\lambda_* \approx 0.1901$ and a guaranteed decrease of the decrement by $\approx 0.2265$. This is a roughly 50\% better performance guarantee than the one derived from \cite[Theorem 5.2.2.1]{Nesterov18book}, and it shall serve as the benchmark against which the short-step method for projectively self-concordant barriers is compared.

%

\medskip

We shall now pass to the analysis of a similar scheme on a projectively self-concordant barrier $f$ with parameter $\gamma = \frac{\nu-2}{\sqrt{\nu-1}}$. Set $F = \nu f$, then $F$ is affinely self-concordant with parameter $\nu$. Let further $p = -\frac{f'}{1 + \langle f',x \rangle}$ be the dual variable. Then $p$ is proportional to $-F'$, and on the central path the dual variable is a positive multiple of $c$. We shall parameterize the central path by the same parameter $\tau$ as in \eqref{affine_central_path}, in order to be able to compare both methods, and denote the corresponding point by $x_*(\tau)$. Besides, we parameterize the path also by a parameter $\kappa$, such that the corresponding point, denoted by $x^*(\kappa)$, obeys the relation $p(x^*(\kappa)) = \kappa c$. Since the method does not compute neither $\kappa$ nor $\tau$, the exact dependence between them is not of interest. Our method alternates Newton-like steps towards a target point on the central path and shifts of the target point along the path corresponding to increments of the parameter $\tau$, or equivalently $\kappa$.

Let us first compute the length of the velocity vector $v$ on the central path in the affine metric $G = f'' - f'f^{\prime T}$. For comparison with the affinely self-concordant case we multiply the affine metric by $\nu$. In view of \eqref{affine_central_path},\eqref{affine_length} we have
\begin{align*}
\sqrt{\nu G(x_*(\tau))[v,v]} =& \sqrt{c^T(F'')^{-1}(F'' - \nu^{-1}F'F^{\prime T})(F'')^{-1}c} \\ =& \sqrt{c^T(F'')^{-1}(F'' - \nu^{-1}\tau^2cc^T)(F'')^{-1}c} = \sqrt{c^T(F'')^{-1}c - \nu^{-1}\tau^2(c(F'')^{-1}c)^2} \\ =& \sqrt{\mu\nu\tau^{-2} - \mu^2\nu\tau^{-2}} = \sqrt{\mu(1-\mu)\nu}\tau^{-1},
\end{align*}
which represents a gain in the length element of a factor $\sqrt{1 - \mu}$ with respect to \eqref{affine_length}.

Let us now pass to the Newton-like step. The quadratic Taylor approximation of $f$ around a point $\hat x \in X^o$ is not projectively equivariant, and we shall replace it by the function
\[ q_{\hat x}(x) = f(\hat x) - \frac12\log\left( (1 - \langle f'(\hat x),x-\hat x \rangle)^2 - G(\hat x)[x-\hat x,x-\hat x] \right),
\]
which shares the function value and the first two derivatives at $\hat x$ with $f$ and is itself projectively self-concordant with the lowest possible parameter value $\gamma = 0$. This introduces a number of differences with respect to the affine case.

Firstly, the domain of definition of the quadratic approximation $q_{\hat x}$ is not $\mathbb R^n$ but the regular convex set
\begin{equation} \label{set_approximation}
\left\{ x \in \mathbb R^n \,\left|\, 1 - \langle f'(\hat x),x-\hat x \rangle > \sqrt{G(\hat x)[x-\hat x,x-\hat x]} \right. \right\},
\end{equation}
which is delimited by a quadric and also serves as an approximation of the set $X$. Secondly, the affine metric of this approximation is not Euclidean but hyperbolic. We therefore have to work with the distance $d_{\hat x}$ it generates rather than a local norm. 
It is explicitly given by
\[ d_{\hat x}(\hat x,x) = \artanh\frac{\sqrt{G[x-\hat x,x-\hat x]}}{1 - \langle f'(\hat x),x-\hat x \rangle}.
\]

We define the projective analogs of the decrement and the Newton step as follows. Let $x_i \in X^o$ be an initial point, and let $\kappa$ be a target parameter. Since we cannot directly access the corresponding point $x^*(\kappa)$, we approximate it by the point $x_f$ obeying the relation $p_{x_i}(x_f) = -\frac{q'_{x_i}(x_f)}{1 + \langle q'_{x_i}(x_f),x_f \rangle} = \kappa c$. The step is then defined by passing from the initial point $x_i$ to the next point $x_f$. The decrement
\[ \rho^{\kappa}(x_i) = d_{x_i}(x_i,x_f)
\]
is defined as the step length, measured in the distance defined by the affine metric of $q_{x_i}$.

In the framework of the short-step method the target parameter is chosen as the maximal value $\kappa$ such that $\rho^{\kappa}(x_i)$ equals some fixed value $\overline{\lambda}$. We shall not compute $\kappa$ explicitly, but rather give a formula for the final point $x_f$. Set $\eta = \tanh^2\overline{\lambda}$, $g = f'(x_i)$, $G = G(x_i)$, and let $\upsilon$ be the maximal root of the quadratic equation
\[ \left( (1-g^TG^{-1}g)^2c^TG^{-1}c + (2-\eta-g^TG^{-1}g)(c^TG^{-1}g)^2 \right)\upsilon^2 + 2(1-\eta)c^TG^{-1}g \cdot \upsilon + g^TG^{-1}g - \eta = 0.
\]
Then the step is given by
\[ x_f - x_i = -\frac{G^{-1}g}{1-g^TG^{-1}g}-\upsilon\left( G^{-1} + \frac{G^{-1}gg^TG^{-1}}{1-g^TG^{-1}g} \right)c.
\]

The analysis of the performance of the step is similar to the affine case. We have to compute an upper bound $\underline{\lambda}$ on the value $\rho^{\kappa}(x_f)$ of the decrement at the final point. Denote the final point corresponding to the target value $\kappa$ and the initial point $x$ by $\chi^{\kappa}(x)$. Let $x(t)$ move along the line segment $\sigma$ linking $x_i$ to $x_f$ with unit velocity until it reaches $x_f$ for some final value $t = T$. For $t > 0$ the point $x_m(t) = \chi^{\kappa}(x(t))$ deviates from $x_f$, and the decrement at the final point equals $\rho^{\kappa}(x_f) = d_{x_f}(x_f,x_m(T))$. Let us track the evolution of the distances $\Delta_f$, $\Delta_m$ from $x(t)$ to $x_f$, $x_m(t)$, respectively, and the angle $\varphi$ between the directions to these points, measured by virtue of the affine metric of $q_{x(t)}$. 
Projective invariance allows to obtain the controlled dynamical system
\[ \frac{d\Delta_f}{dt} = - 1 + \gamma \frac{\sinh2\Delta_f}{2} e_1^TUe_1, \quad \frac{d\Delta_m}{dt} = - \cos\varphi - \gamma \frac{\sinh2\Delta_m}{2} u_{\varphi}^TUu_{\varphi},
\]
\[ \frac{d\varphi}{dt} = \frac{\sin\varphi}{\tanh\Delta_m} + \gamma \frac{\cos\varphi}{\sin\varphi} \cdot \left( e_1^TUe_1 - u_{\varphi}^TUu_{\varphi} \right),
\]
where again $U \in {\cal U}$ is the control. The upper bound $\underline{\lambda}$ is given by the maximal value of $\Delta_m(T)$ over the trajectories of the system with boundary conditions $\Delta_f(0) = \Delta_m(0) = \overline{\lambda}$, $\varphi(0) = 0$, $\Delta_f(T) = 0$.

In general the solution of the problem cannot be presented in closed form, and we consider the 1-dimensional case. Then the system simplifies to
\[ \frac{d\Delta_f}{dt} = - 1 + \gamma \frac{\sinh2\Delta_f}{2} u, \quad \frac{d\Delta_m}{dt} = - 1 - \gamma \frac{\sinh2\Delta_m}{2} u, \quad u \in [-1,1],
\]
where we again admit negative values for $\Delta_m$ and maximize the absolute value $|\Delta_m(T)|$. The solution of the problem is similar to the affine case. The optimal control is $u \equiv -1$ up to the point where $x(t)$ has equal distance to both $x_f$ and $x_m(t)$, and switches to $u \equiv 1$ thereafter. The maximal value of $|\Delta_m(T)|$ as a function of $\overline{\lambda}$ is given by
\[ \underline{\lambda} = \artanh\frac{4\gamma - \gamma^3\tanh^2\overline{\lambda} - 4\gamma\sqrt{(1 - \tanh^2\overline{\lambda})^2 - \gamma^2\tanh^2\overline{\lambda}}}{(\gamma^2 + 4)(1 - \tanh^2\overline{\lambda})}.
\]

In order to obtain the optimal decrement we have to maximize the difference $\overline{\lambda} - \underline{\lambda}$ with respect to $\overline{\lambda}$. The maximizer $\lambda^*$ is given by the roots of a polynomial in $\gamma$ and $\tanh\overline{\lambda}$. We have $\lambda^* \approx -\frac12\log\gamma$ as $\gamma \to 0$ and $\lambda^* \approx \alpha\gamma^{-1}$ as $\gamma \to +\infty$, where $\alpha \approx 0.4166$ is the value of $\lambda^*$ for the affine case. The corresponding lower bound $\lambda_*$ tends to $\artanh(\sqrt{2}-1)$ as $\gamma \to 0$ and behaves like $\lambda^* \approx \beta\gamma^{-1}$ as $\gamma \to +\infty$, where $\beta \approx 0.1901$ is the value of $\lambda_*$ for the affine case.

\begin{figure}
\centering
\includegraphics[width=15.56cm,height=8.71cm]{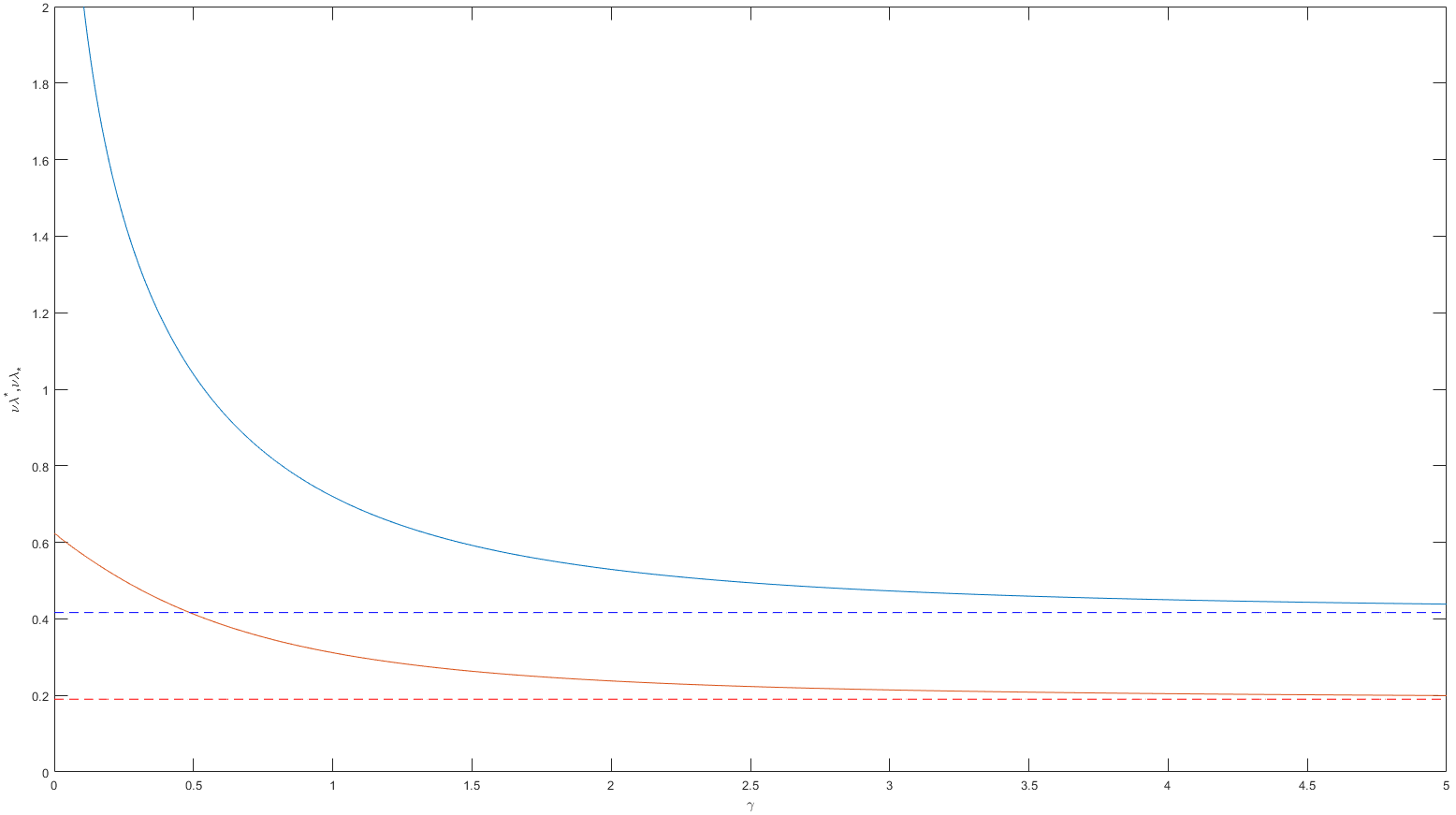}
\caption{Optimal normalized decrements at the initial and final point of a Newton-like step for a projectively self-concordant barrier in one dimension (solid) and the corresponding values for an affinely self-concordant barrier (dashed). }
\label{fig:upper_lower_distances}
\end{figure}

As was mentioned above, for comparison with the affinely self-concordant case we have to normalize the decrements by multiplication by a factor $\sqrt{\nu}$. The normalized values are depicted in Fig.~\ref{fig:upper_lower_distances} as functions of $\gamma$. Since $\gamma \approx \sqrt{\nu}$ for large $\gamma$, the normalized values for a projectively self-concordant barrier tend to the corresponding values for an affinely self-concordant barrier. This property is preserved also for general dimension $n$, as the (normalized) dynamical system for the projective case tends to the dynamical system for the affine case as $\gamma \to +\infty$.

\medskip

Let us compare the performance of the short-step method in dependence on whether the used barrier is considered as projectively or affinely self-concordant. At each iteration the method increases the parameter $\tau$ of the central path in a way such that the length of the jumped over central path segment approximately equals the decrease in the value of the normalized decrement which is achieved by the Newton step. The method makes larger steps if projective self-concordance is taken into account, and there are two different sources of performance gain.

One source is that the same segment of the central path is shorter in the normalized affine metric $\nu G$ by a factor $\sqrt{1-\mu}$, where $\mu(\tau) = \nu^{-1}||F'(x_*(\tau))||^2_{x_*(\tau)}$ indicates the degree of tightness of condition (4) in Definition \ref{def:sc_barrier} on the central path. 

The other source is the ability of the step to decrease the value of the decrement. While for $\gamma \to +\infty$ the performances of the two approaches are asymptotically equal, the decrease with the projective approach tends to infinity for $\gamma \to 0$. This phenomenon can be explained as follows. If $\gamma = 0$, then the barrier coincides with its quadratic approximation, and the solution of the initial optimization problem can be found analytically in a single step, taking us instantly to the infinitely distant limit of the central path.

Let us comment on another observation. In the affinely self-concordant case the analysis of the decrease in the decrement does not use the barrier parameter at all, and therefore implicitly assumes an infinitely large parameter. The reason is that condition (4) in Definition \ref{def:sc_barrier} is not invariant with respect to addition of a linear term to $F$. Projective self-concordance does not tolerate addition of linear terms either, but the larger group of symmetries, namely projective transformations instead of affine ones, allows to set the gradient $f'$ to zero at any given point of the domain by an appropriate choice of coordinates, thus making the theory effectively independent of the magnitude of the gradient.

\section{Outlook} \label{sec:outlook}

In this contribution we presented a new class of barrier functions for convex optimization with a modified self-concordance property, which we called projective self-concordance. It has superior theoretical properties in comparison to the class of classical self-concordant barrier functions. In particular, the Dikin sets are larger and there exists also a quadratic outer approximation of the underlying convex set centered on an arbitrary interior point which can be constructed from the derivatives of the barrier at this point. Moreover, projectively self-concordant barriers admit not only a quadratic approximation of the function, but also furnish a quadratic approximation \eqref{set_approximation} of the domain around an arbitrary interior point. This additional feature might also open possibilities for new classes of interior-point methods. A preliminary analysis of a short-step method shows that taking into account projective self-concordance allows to make larger steps than affine self-concordance alone, for the same barrier.


Another application of projectively self-concordant barriers is the theoretical study of logarithmically homogeneous barriers on cones. A projectively self-concordant barrier on a domain is equivalent to a logaritmically homogeneous barrier on the cone over the domain, but its domain of definition has a smaller dimension. This significantly facilitates the study of 3-dimensional cones, or if properties of the barrier are considered which are essentially defined on two-dimensional sections.

As an example, let us consider the problem of bounding the parameter of convex combinations of self-concordant barriers. More precisely, let $K$ be a regular convex cone, and let $F_1,\dots,F_m$ be logarithmically homogeneous self-concordant barriers on $K$ with parameter $\nu$. A convex combination $F = \sum_{i=1}^m \lambda_iF_i$ of these barriers will in general not be self-concordant. What is then the minimal constant $\vartheta$ such that $\vartheta F$ is a self-concordant barrier? Theorem 5.1.1 in \cite{Nesterov18book} suggests that there is no bound on $\vartheta$ which is uniform over all convex combinations.

However, the problem can be solved affirmatively and exactly when using projective self-concordance. Let $f = \nu^{-1}F|_{X^o}$ be the projectively self-concordant barrier on a compact affine section $X$ of $K$ generated by $F$. Choose an interior point $x \in X^o$ and a tangent vector $h \in T_xX^o$. Then condition (3) in Definition \ref{def:psc_barrier} and the bounds on the derivatives of $f$ in Corollary \ref{cor:bound_sigma} can be written as a set of non-convex constraints on the triple $w = (f'(x)[h],f''(x)[h,h],f'''(x)[h,h,h])$, delimiting a family of nested compact non-convex subsets $B_{\gamma} \subset \mathbb R^3$, parameterized by the self-concordance parameter $\gamma$. Let $w_i$ be the triples defined by barriers $f_i$ with parameter $\gamma$. The triple $w$ defined by a convex combination $f$ of the $f_i$ is an element of the convex hull of $B_{\gamma}$. If we find a number $\gamma'$ such that $B_{\gamma'}$ contains the convex hull of $B_{\gamma}$, then $f$ satisfies condition (3) for the parameter value $\gamma'$. Note that $B_{\gamma}$ depends also on $x,h$ by virtue of the quantities $\sigma_x(h),\pi_x(h)$. However, by projective equivariance of the conditions defining $B_{\gamma}$ the sets are mutually projectively equivalent for different $x,h$, because by an appropriate projective transformation of $X$ and normalization of $h$ we may always achieve $\sigma_x(h) = \pi_x(h) = 1$. But convex hulls are preserved by projective transformations (as opposed to individual convex combinations), and therefore we deal essentially with a single family of nested compact sets $B_{\gamma}$. The dependence of $\gamma' = \frac{\nu'-2}{\sqrt{\nu'-1}}$ on $\gamma = \frac{\nu-2}{\sqrt{\nu-1}}$ can be explicitly computed and translates to a similar dependence of $\nu'$ on $\nu$. A convex combination of logarithmically homogeneous self-concordant barriers on $K$ with parameter $\nu$ is then self-concordant with parameter $\nu'$ after multiplication by $\frac{\nu'}{\nu}$.

\section*{Acknowledgments}

The author would like to thank Prof. Yuri Nesterov for insightful discussions while visiting him in Louvain-la-Neuve in March 2019. He would also like to thank the two anonymous reviewers for helpful comments which led to substantial improvements of the paper.

\bibliographystyle{plain}
\bibliography{interior_point,convexity,affine_geometry}

\end{document}